\pdfoutput=1
\documentclass[a4paper,11pt]{article}
\usepackage[english]{babel}
\usepackage[sc,osf]{mathpazo}
\usepackage[T1]{fontenc}
\usepackage[utf8]{inputenc}
\usepackage{amssymb}
\usepackage{amsthm}
\usepackage{thmtools}
\usepackage{graphicx}
\usepackage{amsmath}
\usepackage[shortlabels]{enumitem}
\usepackage{pifont}
\usepackage{mathtools}
\usepackage{tensor}
\usepackage{stmaryrd}
\usepackage{placeins}
\usepackage{verbatim}
\usepackage{tikz-cd}
\usepackage{tikz}
\usetikzlibrary{shapes.geometric,patterns,arrows,decorations.pathreplacing}
\usepackage[textwidth=1in, textsize = footnotesize]{todonotes}
\usepackage{caption}
\usepackage{sidecap}
\usepackage{subcaption}
\usepackage{accents}
\usepackage{color}
\usepackage{hyperref}
\usepackage{imakeidx}
\usepackage{csquotes}
\usepackage{subfiles}
\usepackage[backend=biber,style=alphabetic,sorting=nty,giveninits]{biblatex}
\renewbibmacro{in:}{%
  \ifentrytype{article}{}{\printtext{\bibstring{in}\intitlepunct}}}
\tikzset{%
    symbol/.style={%
        draw=none,
        every to/.append style={%
            edge node={node [sloped, allow upside down, auto=false]{$#1$}}}
    }
}
\tikzset{
	no line/.style={draw=none,
		commutative diagrams/every label/.append style={/tikz/auto=false}},
	from/.style args={#1 to #2}{to path={(#1)--(#2)\tikztonodes}}}

\tikzset{
    rot90/.style={anchor=south, rotate=90, inner sep=.5mm}
}
\tikzset{
    rot320/.style={anchor=south, rotate=320, inner sep=.5mm}
}
\tikzset{
    rot20/.style={anchor=south, rotate=20, inner sep=.5mm}
}

\addto\extrasenglish{ 
 
 }

\declaretheorem[name=Theorem, numberwithin=section]{theorem}
\declaretheorem[name=Theorem, numbered=no]{theorem*}
\declaretheorem[name=Theorem]{theoremA}

\declaretheorem[name=Lemma, sibling=theorem]{lemma}
\declaretheorem[name=Proposition, sibling=theorem]{proposition}
\declaretheorem[name=Corollary, sibling=theorem]{corollary}

\declaretheorem[style=definition, name=Definition, sibling=theorem]{definition}
\declaretheorem[style=definition, name=Remark, sibling=theorem]{remark}
\declaretheorem[style=definition, name=Example, sibling=theorem]{example}

\DeclareMathOperator{\Fun}{Fun}
\DeclareMathOperator*{\colim}{colim}
\DeclareMathOperator{\id}{id}

\DeclareMathOperator{\Pro}{Pro}

\DeclareMathOperator{\Hom}{Hom}
\DeclareMathOperator{\Fill}{Fill}
\DeclareMathOperator{\Spine}{Sp}
\DeclareMathOperator{\CExp}{Cexp}
\DeclareMathOperator{\sHom}{\mathbf{hom}}
\DeclareMathOperator{\Map}{Map}
\DeclareMathOperator{\Aut}{Aut}
\DeclareMathOperator{\map}{map}

\DeclareMathOperator{\im}{im}

\DeclareMathOperator{\cosk}{cosk}
\DeclareMathOperator{\sk}{sk}

\DeclareMathOperator{\cl}{cl}
\newcommand{\clboun}{\partial_{\mathrm{cl}}}
\newcommand{\wh}{\widehat}
\newcommand{\wt}{\widetilde}
\newcommand{\ol}{\overline}
\newcommand{\ul}{\underline}

\newcommand{\bbN}{\mathbb N}

\newcommand{\bfC}{\mathbf C}
\newcommand{\bfD}{\mathbf D}
\newcommand{\bfE}{\mathcal E}

\newcommand{\bftwo}{\mathbf 2}

\newcommand{\s}{\mathbf s}
\newcommand{\dd}{\mathbf d}
\newcommand{\od}{\mathbf{od}}
\newcommand{\cd}{\mathbf{cd}}

\newcommand{\Fin}{\mathbf{Fin}}
\newcommand{\Set}{\mathbf{Set}}

\newcommand{\Grp}{\mathbf{Grp}}

\newcommand{\Stone}{\mathbf{Stone}}
\newcommand{\Top}{\mathbf{Top}}
\newcommand{\Op}{\mathbf{Op}}
\newcommand{\FinsSet}{\mathbf{sSet}_\mathrm{fin}}

\newcommand{\LdSet}{\mathbf{dL}}

\newcommand{\Ac}[2]{#1 \mhyphen #2}
\newcommand{\AcG}[1]{\Ac{G}{#1}}
\newcommand{\GSet}{\AcG{\Set}}
\newcommand{\GwhSet}{\AcG{\wh\Set}}
\newcommand{\Omegao}{\Omega_\mathrm{o}}
\newcommand{\Omegacl}{\Omega_\mathrm{cl}}

\makeatletter
\newcommand{\oset}[3][0ex]{%
  \mathrel{\mathop{#3}\limits^{
    \vbox to#1{\kern-2\ex@
    \hbox{$\scriptstyle#2$}\vss}}}}
\makeatother

\newcommand{\cofarrow}{\rightarrowtail}
\newcommand{\trivcofarrow}{\oset[-.5ex]{\sim}{\rightarrowtail}}
\newcommand{\wearrow}{\oset[-.16ex]{\sim}{\longrightarrow}}
\newcommand{\fibarrow}{\twoheadrightarrow}
\newcommand{\trivfibarrow}{\oset[-.5ex]{\sim}{\twoheadrightarrow}}

\newcommand{\rlp}{right lifting property }
\newcommand{\llp}{left lifting property }
\newcommand{\wrt}{with respect to }

\renewcommand{\subset}{\subseteq}

\mathchardef\mhyphen="2D

\title{Profinite $\infty$-operads}
\author{Thomas Blom and Ieke Moerdijk}
\date{\today}

\begin{document}

\maketitle
\begin{abstract}
    We show that a profinite completion functor for (simplicial or topological) operads with good homotopical properties can be constructed as a left Quillen functor from an appropriate model category of $\infty$-operads to a certain model category of profinite $\infty$-operads. The construction is based on a notion of lean $\infty$-operad, and we characterize those $\infty$-operads weakly equivalent to lean ones in terms of homotopical finiteness properties. Several variants of the construction are also discussed, such as the cases of unital (or closed) $\infty$-operads and of $\infty$-categories.
\end{abstract}

\section{Introduction}

The goal of this paper is to construct a profinite completion functor for (simplicial or topological) operads, and study its homotopy theoretical properties. Such profinite completions of operads occur, for example, in the work of Horel and of Boavida-Horel-Robertson in their characterization of the Grothendieck-Teichmüller group as the group of self-equivalences of the profinite completion of the little 2-cubes operad \cite{Horel2017ProfiniteOperads,BoavidaHorelRobertson2019Operads}. These self-equivalences are defined in a somewhat ad hoc fashion, and one of our goals is to show that there is a well-defined homotopy theory in the sense of Quillen underlying these. A naive approach to profinite completion immediately stumbles on the problem that profinite completion of spaces does not preserves products, so taking the profinite completion of the spaces of operations in an operad does not produce an operad. Since profinite completion of spaces does sometimes preserve products up to weak equivalence (such as for the spaces of operations of the little 2-cubes operad), the authors of the papers cited above work with $\infty$-operads.

In this paper, we will use the category $\dd\Set$ of dendroidal sets, equipped with the so-called operadic model structure of \cite{CisinskiMoerdijk2011DendroidalSets}, as a model for $\infty$-operads. This model category is Quillen equivalent to the model categories of simplicial or topological operads. A natural candidate for a category to model the homotopy theory of profinite $\infty$-operads would then be the category of dendroidal objects valued in profinite sets, or equivalently, dendroidal Stone spaces. This category can be shown to be equivalent to the pro-category of the full subcategory of $\dd\Set$ spanned by the so-called \emph{lean} dendroidal sets. A dendroidal set is called lean if it is the $n$-coskeleton of a degreewise finite dendroidal set. A lean dendroidal set that is furthermore fibrant in the operadic model structure will be called a \emph{lean $\infty$-operad}. The forgetful functor sending a dendroidal profinite set to its underlying dendroidal set admits a left adjoint, which will be called \emph{profinite completion}. Our main result is the following.

\begin{theoremA}\label{TheoremA}
	There exists a fibrantly generated model structure on the category of dendroidal profinite sets such that the weak equivalences between lean $\infty$-operads are exactly those of the operadic model structure on dendroidal sets, and for which the profinite completion functor is a left Quillen functor from the operadic model structure to this new model structure.
\end{theoremA}

This model structure will be called the \emph{model category of profinite $\infty$-operads}. The generating (trivial) fibrations of this model structure are the (trivial) fibrations between lean $\infty$-operads in the operadic model structure on dendroidal sets. We refer to \autoref{theorem:ProfiniteOperadicModelStructure} for a more precise formulation of this theorem and some variations, namely versions of this result for open and closed dendroidal sets and for other model structures than the operadic one (cf. \autoref{remark:CovariantPicard}). In particular, we will obtain as a consequence a similar model structure for profinite $\infty$-categories, related to the Joyal model structure for $\infty$-categories by a Quillen adjunction in which the left adjoint is the profinite completion functor.

Key to understanding the model category of profinite $\infty$-operads is a good grasp of which $\infty$-operads are (equivalent to) lean $\infty$-operads. We will give a precise homotopical criterion for being equivalent to a lean $\infty$-operad in \autoref{ssec:LeanInfinityOperads}. Namely, we show that an $\infty$-operad is equivalent to a lean $\infty$-operad precisely when it has finitely many colours up to homotopy and its spaces of operations are all $\pi$-finite and contractible in sufficiently high arity (see \autoref{def:PiFiniteInftyOperad} for a precise formulation). Such $\infty$-operads will be called \textit{$\pi$-finite}.

Our strategy is similar to the one in our earlier paper \cite{BlomMoerdijk2020SimplicialProV1}, but the situation here is more involved for two reasons. The first of these is that strictly speaking, the category of dendroidal sets is not simplicially enriched. This does not turn out to be a real issue, since there is still a canonical way to define a simplicial hom satisfying all properties we will need (cf. \autoref{ssec:dSetsAndSpaces}). The second reason is that not all objects are cofibrant in the operadic model structure, but only the so-called \emph{normal} ones. This will make it harder to give a correct definition of the weak equivalences between dendroidal profinite sets and in particular forces us to study the normal monomorphisms of dendroidal profinite sets in detail.

Our method for constructing the model structure is surprisingly general and flexible. In a sequel to this paper \cite{BlomMoerdijk2020ProfiniteInftyOperadsPartIITA}, we will show that there is a similar profinite version of the complete (in the sense of Rezk) Segal model structure on dendroidal \emph{spaces}, Quillen equivalent to the model structure for profinite $\infty$-operads constructed in this paper. This complete Segal type model structure will have some convenient properties that will give us a Dwyer-Kan style characterization of the weak equivalences between pofinite $\infty$-operads as the essentially surjective and fully faithful maps, appropriately defined. This model structure will also enable us to obtain the following description of the $\infty$-category associated to the model category of profinite $\infty$-operads, which we already wish to single out here since it makes precise the role of the $\pi$-finite $\infty$-operads within the homotopy theory of all profinite $\infty$-operads. Let $\Pro$ denote the pro-category of an $\infty$-category defined as in the dual of \cite[Definition 5.3.5.1]{Lurie2009HTT} and let $\Op_\infty$ denote the $\infty$-category of $\infty$-operads.

\begin{theorem*}[\cite{BlomMoerdijk2020ProfiniteInftyOperadsPartIITA}]
	The underlying $\infty$-category of the model category of profinite $\infty$-operads is equivalent to $\Pro(\Op^\pi_\infty)$, where $\Op^\pi_\infty \subset \Op_\infty$ denotes the full sub-$\infty$-category spanned by the $\pi$-finite $\infty$-operads.
\end{theorem*}

The results of this paper have been written in such a way that they also apply to open and closed dendroidal sets. The reason for including closed dendroidal sets is that the model structure for profinite $\infty$-operads mentioned above does not appear to be satisfactory when working with unital operads. This has to do with there not existing ``enough'' unital lean $\infty$-operads. For example, we will see in \autoref{remark:AssociativeOperadNotProfinite} below that the (unital) associative operad can't be written as an inverse limit of lean $\infty$-operads and in particular that it is not a a profinite $\infty$-operad. This issue can be resolved by working with closed dendroidal sets instead. The category of closed dendroidal sets admits a model structure Quillen equivalent to that of unital simplicial operads, and a model structure on closed dendroidal profinite sets analogous to that of \autoref{TheoremA} is constructed in \autoref{theorem:ProfiniteOperadicModelStructure} below. The condition of being $n$-coskeletal for closed dendroidal sets is different in nature than for general dendroidal sets. In \autoref{ssec:LeanInfinityOperads}, it will be shown that this results in a different characterization of those closed $\infty$-operads which are equivalent to lean closed $\infty$-operads (i.e. to fibrant closed dendroidal sets that are coskeletal and degreewise finite). This characterization is more natural when working with unital operads, as is illustrated by the fact that the unital associative operad is a lean closed $\infty$-operad (see \autoref{exx:AssociativeOperadCoskeletal} below). In particular, working with closed dendroidal sets solves the issue of there ``not being enough'' unital lean $\infty$-operads.

\paragraph{Overview of the paper.}

In \autoref{sec:DendroidalSets}, we discuss some generalities on dendroidal sets and prove our homotopy-theoretic characterization of lean $\infty$-operads. We next recall some basic facts on pro-categories and prove some elementary results on profinite sets, or equivalently Stone spaces, in \autoref{sec:ProCategories}. These elementary results on profinite sets are then used to study normal monomorphisms of dendroidal profinite sets in \autoref{sec:NormalMonos}. Finally, in \autoref{sec:TheModelStructure} we construct the model structure for profinite $\infty$-operads and its open and closed counterparts, and we deduce some basic properties.

\section{Dendroidal sets} \label{sec:DendroidalSets}

We start this section by recalling some basic definitions and facts on dendroidal sets. We will then introduce skeleta and coskeleta of dendroidal sets, which will be used throughout many proofs in this paper. Particularly important in the construction of the model structure for profinite $\infty$-operads in \autoref{sec:TheModelStructure} are the so-called \emph{lean} dendroidal sets; that is, the degreewise finite and coskeletal dendroidal sets. For this reason, we conclude this section by giving a precise characterization of which $\infty$-operads are weakly equivalent to lean $\infty$-operads.

\subsection{Preliminaries on dendroidal sets} \label{ssec:dSetsAndSpaces}

For the convenience of the reader, we will quickly review the basic definitions and facts concerning dendroidal sets. For details, we refer to \cite{HeutsMoerdijk2020Trees} and the references cited there.

\paragraph{Trees.} The theory is based on a category $\Omega$ of trees. Its objects are finite rooted trees, where the top edges can be ``open'', in which case they are called \emph{leaves}, or ``closed'', in which case they are called \emph{stumps}. Here is a typical example of such a tree.
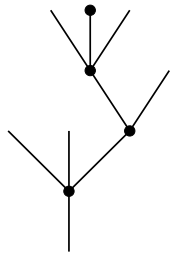
\[\scalebox{0.8}{\begin{tikzpicture} 
    \draw[fill] (0,1) circle [radius=0.085];
    \draw[fill] (1,2) circle [radius=0.085];
    \draw[fill] (0.35,3) circle [radius=0.085];
    \draw[fill] (0.35,4) circle [radius=0.085];
    
    \draw[thick] (0,0) -- (0,1);
    \draw[thick] (0,1) -- (-1,2);
    \draw[thick] (0,1) -- (0,2);
    \draw[thick] (0,1) -- (1,2);
    \draw[thick] (1,2) -- (1.65,3);
    \draw[thick] (1,2) -- (0.35,3);
    \draw[thick] (0.35,3) -- (-0.30,4);
    \draw[thick] (0.35,3) -- (0.35,4);
    \draw[thick] (0.35,3) -- (1,4);
\end{tikzpicture}}\]

The edges of the tree other than the root and the leaves are called \emph{inner edges}. Each such tree $T$ generates a coloured symmetric operad $\Omega(T)$ whose colours are the edges of $T$, and whose operations are generated by the vertices of $T$. The morphisms $S \to T$ in the category $\Omega$ are the operad maps $\Omega(S) \to \Omega(T)$. These can be described as compositions of ``(elementary) face maps'', of isomorphisms, and of ``(elementary) degeneracies''; see \cite[\S 3.3]{HeutsMoerdijk2020Trees}. The simplex category $\Delta$ embeds into $\Omega$ as the full subcategory on the linear trees with a single leaf, which are the trees of the form depicted below.
\[\scalebox{0.8}{\begin{tikzpicture} 
    \draw[thick] (0,0) -- (0,1.25);
    \draw[thick] (0,1.85) -- (0,2.6);
    
    \draw[fill] (0,0.4) circle [radius=0.085];
    \draw[fill] (0,0.9) circle [radius=0.085];
    \draw[fill] (0,2.2) circle [radius=0.085];
    
    \node at (0,1.65) {$\vdots$};
\end{tikzpicture}}\]
In fact, writing $\eta$ for the unique linear tree with no vertices, the inclusion $\Delta \hookrightarrow \Omega$ factors as $\Delta \xrightarrow{\cong} \Omega/\eta \to \Omega$, where the first functor is an isomorphism.

A \emph{dendroidal set} is a set-valued presheaf on $\Omega$, and the category of these is denoted $\dd\Set$. For such a presheaf $X$, the value at a tree $T$ will be denoted by $X_T$. We write $\Omega[T]$ for the presheaf represented by $T$ and $\partial \Omega[T]$ for its boundary, i.e. the subobject of $\Omega[T]$ obtained as the union indexed by all proper faces $S \rightarrowtail T$. We will often simply write $T$ for the representable presheaf $\Omega[T]$.

The inclusion $i \colon \Delta \hookrightarrow \Omega$ induces an adjunction
\[i_! : \s\Set \rightleftarrows \dd\Set : i^*. \]
The left adjoint $i_!$ is fully faithful, and its image can be identified with $\dd\Set/\eta$. Here $\eta \in \Omega$ is identified with the presheaf that it represents.

\paragraph{The operadic model structure.} A dendroidal set $X$ is called \emph{normal} if for each tree $T$, the action of $\Aut(T)$ on $X_T$ coming from the presheaf structure is free. Similarly, a monomorphism $X \to Y$ is called \emph{normal} if for each $T$, the action of $\Aut(T)$ on the complement of the image of $X_T \to Y_T$ is free. The category $\dd\Set$ carries the so-called \emph{operadic model structure} \cite{CisinskiMoerdijk2011DendroidalSets}, whose cofibrations are these normal monomorphisms and whose fibrant objects are the \emph{$\infty$-operads}, i.e. those dendroidal sets having the right lifting property with respect to inner horn inclusions. These are inclusions of the form $\Lambda^e[T] \cofarrow \Omega[T]$, where $e$ is an inner edge of $T$ and $\Lambda^e[T]$ is the union of all the proper faces $S \cofarrow T$ that contain the edge $e$. Under the identification $\s\Set = \dd\Set/\eta$, this operadic model structure restricts to the Joyal model structure on $\s\Set$ describing the homotopy theory of \emph{$\infty$-categories}.

The interest of the operadic model structure lies in the fact that it is Quillen equivalent to the natural model structure on the category of simplicial (coloured) operads. This Quillen equivalence is given by a pair of adjoint functors
\begin{equation}\label{eq:CoherentNerveQuillenEquiv}
	w_! : \dd\Set \rightleftarrows \s\Op : w^*,
\end{equation} 
completely determined up to natural isomorphism by the effect of $w_!$ on representables, where it is defined by sending $\Omega[T]$ to the Boardman-Vogt resolution $W \Omega(T) = w_!(\Omega[T])$ of the operad $\Omega(T)$ freely generated by $T$. The right adjoint $w^*$ of this Quillen pair is called the \emph{homotopy-coherent nerve}. When restricted to dendroidal sets over $\eta$, this Quillen equivalence recovers the well-known one between the Joyal model structure and simplicial categories.

\paragraph{The tensor product and the simplicial hom.} The category of dendroidal sets admits a tensor product $\otimes$ that is related to the Boardman-Vogt tensor product of operads, see e.g. Chapter 4 of \cite{HeutsMoerdijk2020Trees}. This tensor product is symmetric and admits a right adjoint in both variables, which we will denote as an exponential and call the \emph{internal hom}. For two dendroidal sets $X$ and $Y$, we define the \emph{simplicial hom} $\sHom(X,Y)$ by
\[\sHom(X,Y)_\bullet = i^*(Y^X) \cong \Hom(X \otimes i_!\Delta[\bullet], Y) \cong \Hom(X,Y^{i_!\Delta[\bullet]}). \]
Strictly speaking, this does not define a simplicial enrichment of the category $\dd\Set$, due to the subtle fact that the relevant tensor product is only associative up to weak equivalence. However, $\sHom(-,-)$ still defines a functor $\dd\Set^{op} \times \dd\Set \to \s\Set$ and this functor interacts well with the operadic model structure on $\dd\Set$; namely, one can show that for any normal monomorphism $X \cofarrow Y$ and any fibration $L \fibarrow K$ in the operadic model structure on $\dd\Set$, the pullback-power map
\[\sHom(Y,L) \to \sHom(X,L) \times_{\sHom(X,K)} \sHom(Y,K) \]
is a fibration in the Joyal model structure on $\s\Set$, which is trivial if either $X \cofarrow Y$ or $L \fibarrow K$ is a weak equivalence in the operadic model structure on $\dd\Set$. This means that, in practice, many of the techniques commonly used in simplicial model categories still apply to $\dd\Set$ (with respect to the Joyal model structure on $\s\Set$).

The simplicial hom can be used to define a model for the mapping space between dendroidal sets. The inclusion of the category of Kan complexes into the category of $\infty$-categories admits a right adjoint $k$ that sends sends an $\infty$-category $X$ to the maximal Kan complex $kX$ that it contains (\cite[Corollary 1.5]{Joyal2002QuasiKan}). For any normal dendroidal set $X$ and any $\infty$-operad $Y$, the Kan complex $\Map(X,Y) = k \sHom(X,Y)$ is a model for the mapping space from $X$ to $Y$. The functor $k$ sends equivalences of $\infty$-categories to homotopy equivalences of Kan complexes and categorical fibrations to Kan fibrations, so the properties of the simplicial hom mentioned above imply that for any normal monomorphism between normal dendroidal sets $X \cofarrow Y$ and any fibration of $\infty$-operads $L \fibarrow K$, the map
\[\Map(Y,L) \to \Map(X,L) \times_{\Map(X,K)} \Map(Y,K) \]
is a Kan fibration that is trivial whenever $X \cofarrow Y$ or $L \fibarrow K$ is a weak equivalence. Throughout this paper, whenever we write $\Map(X,Y)$, we will always refer to this particular model for the mapping space.

\paragraph{Spaces of operations.} For an $\infty$-operad $X$, we call $X_\eta$ its \emph{set of colours}. For a tuple of colours $c_1, \ldots, c_n,d \in X_\eta$, the \emph{space of operations} $X(c_1,\ldots,c_n;d)$ is defined as the pullback
\[\begin{tikzcd}
X(c_1,\ldots,c_n;d) \ar[r] \ar[d] \ar[dr, phantom, very near start, "\lrcorner"] & \sHom(\Omega[C_n],X) \ar[d] \\
\Delta[0] \ar[r,"(c_1 {,} \ldots {,} c_n {,} d)"] & \sHom(\partial \Omega[C_n], X) \cong (i^*X)^{n+1}.
\end{tikzcd}\]
Here $C_n$ denotes the $n$-corolla, the unique tree with one vertex and $n$ leaves. For the isomorphism on the bottom right, note that $\partial \Omega[C_n] \cong \sqcup_{0 \leq j \leq n} \eta$ and that $\sHom(\eta,X) = i^* X$. The map $X(c_1,\ldots,c_n;d) \to \sHom(\Omega[C_n],X)$ lands in $\Map(\Omega[C_n],X) = k\sHom(\Omega[C_n],X)$, hence the space of operations $X(c_1,\ldots,c_n;d)$ agrees with the fiber of $\Map(\Omega[C_n],X) \fibarrow \Map(\partial \Omega[C_n], X)$ above $(c_1,\ldots,c_n,d)$ (cf. Theorem 6.51(b) and Remark 9.43 of \cite{HeutsMoerdijk2020Trees}).

A map $f \colon X \to Y$ of $\infty$-operads is called \emph{essentially surjective} if $\pi_0(\Map(\eta,X)) \to \pi_0(\Map(\eta,Y))$ is surjective, and \emph{fully faithful} if for any tuple $c_1,\ldots,c_n,d \in X_\eta$, the map $X(c_1,\ldots,c_n;d) \to Y(f(c_1),\ldots,f(c_n);f(d))$ is a weak equivalence of Kan complexes. One can show that the weak equivalences between $\infty$-operads are precisely the essentially surjective and fully faithful maps (cf. \cite[Theorem 9.45]{HeutsMoerdijk2020Trees}).

\paragraph{Open dendroidal sets.} A tree is called \emph{open} if it has no stumps. The full subcategory of $\Omega$ spanned by the open trees is denoted $\Omegao$, and the inclusion is denoted $o \colon \Omegao \hookrightarrow \Omega$. We write $\od\Set$ for the category of ($\Set$-valued) presheaves on $\Omegao$, and refer to its objects as \emph{open dendroidal sets}. The functor $o$ induces an adjoint pair
\begin{equation}\label{eq:oshriekAdjunction}
	o_! : \od\Set \rightleftarrows \dd\Set : o^*
\end{equation}
for which $o_!$ is a full embedding and $o^*$ is simply the restriction along $o$. Under the full embedding $o_!$, the category $\od\Set$ is identified with the full subcategory of $\dd\Set$ spanned by those dendroidal sets $X$ that have the property that $X_T = \varnothing$ whenever $T$ is not open. The category $\od\Set$ can also be identified with the slice category $\dd\Set/O$ over the subobject $O \subset *$ of the terminal object whose value $O(T)$ at a tree $T$ is non-empty if and only if $T$ is open. Note that the inclusion $i_! \colon \s\Set \to \dd\Set$ factors through $\od\Set$, and that the tensor product of two open dendroidal sets (viewed as object of $\dd\Set$) is again open. In other words, $o_!$ preserves the tensor product. In particular, tensors of open dendroidal sets with simplicial sets are automatically open and left adjoint to the simplicial hom of $\dd\Set$ restricted to $\od\Set$.

The operadic model structure restricts in the usual way to the slice $\od\Set \simeq \dd\Set/O$, making \eqref{eq:oshriekAdjunction} into a Quillen pair, and the Quillen equivalence given by $w_!$ and $w^*$ restricts to a Quillen equivalence
\[\mathring{w}_! : \od\Set \rightleftarrows \mathbf{o}\s\Op : \mathring{w}^* \]
between open dendroidal sets and ``open'' simplicial operads, i.e. (coloured) simplicial operads without nullary operations.

\paragraph{Closed dendroidal sets.} The following is a quick review of the theory of closed dendroidal sets; for details the reader is referred to \cite{Moerdijk2018Closed}.

A tree is called \emph{closed} if it has no leaves. The full subcategory of $\Omega$ spanned by the closed trees is denoted $\Omegacl$, and the inclusion is denoted $u \colon \Omegacl \hookrightarrow \Omega$. This inclusion has a left adjoint $\cl \colon \Omega \to \Omegacl$ sending a tree $T$ to its \emph{closure} $\ol T = \cl(T)$, constructed by capping each leaf of $T$ by a stump. The category of \emph{closed dendroidal sets}, i.e. presheaves on $\Omegacl$, is denoted $\cd\Set$. We will write $\Omegacl[T]$ for the presheaf represented by the closed tree $T$ and $\clboun \Omegacl[T] \subset \Omegacl[T]$ for its boundary, i.e. the union of all \emph{closed} faces of $T$. The adjoint pair $\cl \dashv u$ induces adjoint functors
\[ \cl_! \dashv u_! = \cl^* \dashv u^* = \cl_* \dashv u_* \]
between $\dd\Set$ and $\cd\Set$. The functor $u_! \colon \cd\Set \to \dd\Set$ is fully faithful and its essential image consists of those dendroidal sets $X$ with the property that for every tree $T$, the map $X_{\ol T} \to X_T$ induced by $T \hookrightarrow \ol T$ is an isomorphism. Identifying $\cd\Set$ with this full subcategory, the tensor product of closed dendroidal sets is again closed and hence restricts to a tensor product on $\cd\Set$. More precisely, this tensor product is defined by $X \otimes Y = \cl_!(u_! X \otimes u_! Y)$. For a simplicial set $M$ and a closed dendroidal set $X$, the tensor $X \otimes M$ is defined as $X \otimes \cl_! i_! M$. A simplicial hom is then defined by the formula
\[\sHom(X,Y)_\bullet = \Hom_{\cd\Set}(X \otimes \Delta[\bullet], Y)\]
for closed dendroidal sets $X$ and $Y$. One can show that $X \otimes \cl_! i_! M$ and $u_! X \otimes i_! M$ agree under the inclusion $u_! \colon \cd\Set \hookrightarrow \dd\Set$, so in particular $\sHom(X,Y)$ agrees with the simplicial hom of $\dd\Set$ under this inclusion.

The category $\cd\Set$ carries a variant of the operadic model structure whose fibrant objects are called \emph{closed $\infty$-operads}. These are the closed dendroidal sets having the \rlp \wrt the ``very inner horns'' $\Lambda^e_{\cl}[T] \cofarrow \Omegacl[T]$, i.e. those inner horns where $e$ is an inner edge of $T$ \emph{not} immediately below a stump. (We have written $\Lambda^e_{\cl}[T]$ to distinguish it from $\Lambda^e[T]$, as $\Lambda^e_{\cl}[T] \subset \clboun \Omegacl[T]$ only involves \emph{closed} faces of $T$ containing the edge $e$.)

This model structure is Quillen equivalent to a (Reedy style) model structure on the category of \emph{closed} or \emph{unital} simplicial operads by a variant of the Quillen equivalence \eqref{eq:CoherentNerveQuillenEquiv}, denoted
\[ \ol w_! : \cd\Set \rightleftarrows \mathbf{u}\s\Op : \ol w^*. \]
The weak equivalences between closed $\infty$-operads can again be characterized as essentially surjective and fully faithful maps, but the definition of the spaces of operations (like $X(c_1,\ldots,c_n;d)$ above) in a closed $\infty$-operad has to be modified slightly, by replacing the corolla $C_n$ by its closed variant $\ol C_n$ and the boundary $\partial \Omega[C_n]$ by the closure $\cl_!(\partial \Omega[C_n]) = \coprod_{i = 0}^n \ol \eta$. (Note that the boundary $\clboun \Omegacl[\ol C_n]$ is much larger than $\coprod \ol \eta$.)

The model structures on $\dd\Set$, $\od\Set$ and $\cd\Set$ are related by several Quillen pairs: we have already mentioned that the embedding $o_! \colon \od\Set \to \dd\Set$ is left Quillen, and the same is true for $u_! = \cl^* \colon \cd\Set \to \dd\Set$. The functor $\cl_! \colon \dd\Set \to \cd\Set$ not left Quillen, but the composition $h = \cl \circ o \colon \Omegao \to \Omegacl$ does induce a Quillen pair. However, now $h^*$ is the left Quillen functor (rather than $h_!$), forming a Quillen pair with its right adjoint $h_* = \cl_* o_* = u^* o_*$.

\subsection{Skeleta and coskeleta of dendroidal sets} \label{ssec:LeanDendroidalSets}

The following discussion of skeleta and coskeleta of dendroidal sets goes through without any change for closed and open dendroidal sets; that is, presheaves on $\Omegacl$ and $\Omegao$, respectively. For simplicity of exposition, we only discuss general dendroidal sets. Given an object $T$ in $\Omega$, we shall write $|T|$ for the sum of the number of non-root edges and the number of vertices in $T$, and call it the \emph{size} of $T$. 

\begin{remark}
	Note that the image of $[n]$ under the embedding $i \colon \Delta \hookrightarrow \Omega$ has size $2n$. This somewhat inconvenient property could be fixed by using a different notion of ``size'', for example by counting only the number of vertices. However, this would imply that the number of trees of size $\leq n$ is infinite for $n \geq 1$. This leads to problems in many of our proofs below. When restricting one's attention the category of closed trees this issue disappears, so a better notion of ``size'' would be to only count the number of vertices in this case. However, to keep our proofs uniform, we use the same definition of size in all three cases, namely the sum of the number of vertices and non-root edges. Our choice of the ``size'' of a tree implies that the functor $\sk_n$ defined below does \textbf{not} agree with the $n$-skeleton of a dendroidal set as defined in \cite[\S 3.6]{HeutsMoerdijk2020Trees}, since the filtration $\cup_n \Omega_{\leq n}$ of $\Omega$ used there is based on counting only the vertices of a tree.
\end{remark}

Let $\Omega_{(n)} \subset \Omega$ denote the full subcategory of trees of size $\leq n$. One easily checks that for each $n$, this full subcategory is finite. The inclusion $\Omega_{(n)} \hookrightarrow \Omega$ defines skeleton and coskeleton functors
\[\sk_n, \cosk_n \colon \dd\Set \to \dd\Set \]
by left and right Kan extension. Explicitly, for a dendroidal set $X$,
\begin{align*}
    \sk_n(X)_T = \colim_S X_S,  \quad &\text{colimit over } T \to S \text{ with }|S| \leq n,\\
    \cosk_n(X)_T = \lim_S X_S, \quad &\text{limit over } S \to T \text{ with } |S| \leq n.
\end{align*}
Since $\Omega_{(n)}$ is finite, these are finite (co)limits. Note that for a given $n \in \bbN$ and tree $T$, these colimits and limits can be taken over the (co)final subcategories of degeneracies $T \twoheadrightarrow S$ with $|S| \leq n$ and of faces $S \rightarrowtail T$ with $|S| \leq n$, respectively.

It immediately follows that for any tree $T$ of size $n$, there is a canonical isomorphism $\sk_{n-1} \Omega[T] \cong \partial \Omega[T]$, where we recall that $\partial \Omega[T]$ is the subobject of the representable presheaf $\Omega[T]$ that is the union of all proper faces $S \rightarrowtail T$. More generally, the $n$-skeleton of a dendroidal set $X$ admits the following simple description: the counit $\sk_n(X) \to X$ gives an isomorphism of $\sk_n(X)$ onto the subobject of $X$ consisting of those $x \in X_T$ for which there exists a degeneracy $\sigma \colon T \twoheadrightarrow S$ with $|S| \leq n$ and an element $y \in X_S$ such that $\sigma^*(y) = x$. In particular, we will identify $\sk_n(X)$ with this subobject of $X$. This description is a direct consequence from Corollary 6.10 of \cite{BergerMoerdijk2011ExtensionReedy}, noting that the ``size'' $|T|$ of a tree $T$ defined above makes $\Omega$ into an EZ-category in the sense of \cite[Definition 6.7]{BergerMoerdijk2011ExtensionReedy}.

By adjunction, there is a natural correspondence of morphisms of dendroidal sets
\[\Hom(\sk_n X, Y) \cong \Hom(X, \cosk_n Y) \]
for any $n \geq 0$ and any dendroidal sets $X$ and $Y$. This yields an alternative description of $\cosk_n(X)$; namely
\[\cosk_n(X)_T \cong \Hom(\sk_n\Omega[T], X).\]
The inclusions $\Omega_{(n)} \subset \Omega_{(n+1)}$ induce the \emph{skeletal filtration}
\[\sk_0(X) \hookrightarrow \ldots \hookrightarrow \sk_n(X) \hookrightarrow \sk_{n+1}(X) \hookrightarrow \ldots \to \colim_n \sk_n(X) = X \]
and \emph{coskeleton tower}
\[ X = \lim_n \cosk_n(X) \to \ldots \to \cosk_{n+1}(X) \to \cosk_{n}(X) \to \ldots \to \cosk_0(X) \]
of a dendroidal set $X$. By \cite[Proposition 7.3.(iii)]{BergerMoerdijk2011ExtensionReedy}, for any normal dendroidal set $X$ and any $n \geq 0$, the map $\sk_{n-1} X \hookrightarrow \sk_n X$ is obtained as a pushout of a coproduct of boundary inclusions $\partial \Omega[T] \cofarrow \Omega[T]$, where $T$ is of size $n$.

It follows from descriptions of the $n$-skeleton and $n$-coskeleton given above that $\sk_n X\hookrightarrow X$ and $X \to \cosk_n X$ induce isomorphisms $\sk_n(X)_T \cong X_T$ and $X_T \cong \cosk_n(X)_T$ if $|T| \leq n$. We will say that $X \to Y$ is an \emph{isomorphism on $n$-skeleta} if $X_T \to Y_T$ is an isomorphism for any tree $T$ of size $\leq n$.

Recall that a dendroidal set is called \emph{finite} if it has finitely many non-degenerate elements. By the above characterization of $\sk_n X$ as a subobject of $X$, one can rephrase this by saying that a dendroidal set $X$ is finite if and only if it is degreewise finite and there exists an $n$ such that $\sk_n X \to X$ is an isomorphism. The following dual notion will play an important role in this paper.

\begin{definition}\label{def:LeanDendroidalSet}
A dendroidal set $X$ is called \emph{lean} if
\begin{enumerate}[(a)]
    \item $X$ is degreewise finite; i.e. for each tree $T$ the set $X_T$ is finite, and
    \item $X$ is coskeletal; i.e. there exists an $n \geq 0$ for which the map $X \to \cosk_n X$ is an isomorphism.
\end{enumerate}
If $X$ is furthermore fibrant in the operadic model structure, then it is called a \emph{lean $\infty$-operad}.
The full subcategory of degreewise finite dendroidal sets is denoted $\dd\Fin\Set$ and we write $\LdSet$ for the full subcategory spanned by the lean dendroidal sets.
\end{definition}

\begin{remark}\label{remark:AlternativeDefinitionLean}
Since $\Omega$ has finite hom-sets, it follows that for any finite subcategory $\bfC \subset \Omega$ and any functor $\bfC^{op} \to \Fin\Set$, the right Kan extension along $\bfC^{op} \hookrightarrow \Omega^{op}$ exists and is defined in terms of finite limits. Since any finite subcategory is contained in $\Omega_{(n)}$ for some $n$, it follows that a dendroidal set $X$ is lean if and only if there exist a finite subcategory $\bfC \subset \Omega$ and a functor $Z \colon \bfC^{op} \to \Fin\Set$ such that $X$ is isomorphic to the right Kan extension of $Z$ along $\bfC^{op} \hookrightarrow \Omega^{op}$.
\end{remark}

\begin{remark}\label{remark:PushforwardPreservesLean}
	Recall the functors $o$, $\cl$, $u$ and $h$ from \autoref{ssec:dSetsAndSpaces}. One can show that the induced functors $o_*$, $\cl_*$, $u_*$ and $h_*$ preserve lean objects. We show this for $\cl_*$, the argument for the other functors is similar. Let $X$ be a lean dendroidal set, say that $X$ is the right Kan extension of $Y \colon \bfC^{op} \to \Fin\Set$ along $\bfC^{op} \hookrightarrow \Omega^{op}$ for some finite subcategory $\bfC$ of $\Omega$. Let $\bfD \subset \Omegacl$ be a finite subcategory that contains the image of $\bfC$ under $\cl$. Then $\cl_* X$ is obtained by right Kan extending $Y$  along $\bfC^{op} \to \bfD^{op}$ and then along $\bfD^{op} \hookrightarrow \Omegacl^{op}$. Since both $\bfC$ and $\bfD$ are finite, combining this with the previous remark implies that $\cl_* X$ is lean.
\end{remark}

\begin{remark}\label{remark:PullbackOpenPreservesLean}
	Perhaps somewhat surprising, the functors $o_! \colon \od\Set \to \dd\Set$ and $o^* \colon \dd\Set \to \od\Set$ both also preserve lean objects. For $o^*$, this can be deduced from the fact that if $T$ is an open tree and $S \to T$ is any map of trees, then $S$ must be open as well. For $o_!$, note that this functor is given by
	\[		(o_!X)_T = \begin{cases}
			X_T &\quad \text{if }T \text{ is open}\\
			\varnothing &\quad \text{else.}
			\end{cases}
	\]
	One can deduce from this explicit description that if $X$ is $n$-coskeletal with $n \geq 1$, then $o_! X$ is $n$-coskeletal, so $o_!$ preserves lean objects. We leave the details to the reader. 
\end{remark}

A simplicial set $X$ will be called \emph{lean} if it is degreewise finite and there exist an $n \in \bbN$ such that $X \to \cosk_n X$ is an isomorphism. Here $\cosk_n X$ is defined by restricting $X$ to a presheaf on $\Delta_{\leq n}$ and then right Kan extending along $\Delta^{op}_{\leq n} \hookrightarrow \Delta^{op}$.

\begin{remark}\label{remark:UnderlyingInftyCatPreservesLean}
	Recall the inclusions $\Delta \hookrightarrow \Omegao$ and $\Delta \hookrightarrow \Omega$, which we both denote by $i$. As above, one can deduce that both restriction functors $i^* \colon \od\Set \to \s\Set$ and $i^* \colon \dd\Set \to \s\Set$ preserve lean objects. However, unlike in the case of $o_!$, their left adjoints $i_!$ do not preserve lean objects.
\end{remark}

The following observations about exponentials of dendroidal sets will be important throughout this paper. Recall from \autoref{ssec:dSetsAndSpaces} that the tensor product $\otimes$ of dendroidal sets induces an internal hom that does not agree with the cartesian exponential. The cartesian exponential will be denoted by $\CExp(Y,X)$ to distinguish it from this internal hom $X^Y$.

\begin{lemma}\label{lemma:CartesianExponentialLean}
Let $E$ be a degreewise finite dendroidal set and let $X$ be a lean dendroidal set. Then the cartesian exponential $\CExp(E,X)$ is again lean.
\end{lemma}

\begin{proof}
Let $n$ be such that $X$ is $n$-coskeletal. To see that $\CExp(E,X)$ is degreewise finite, note that
\[\CExp(E,X)_T \cong \Hom(E \times \Omega[T], \cosk_n(X))\]
where the right-hand side is finite since $\cosk_n(X)$ is the right Kan extension of a functor $\Omega_{(n)} \to \Fin\Set$, where $\Omega_{(n)}$ is a finite category. To see that $\CExp(E,X)$ is coskeletal, let $Y$ be any dendroidal set. Then
\begin{align*}
    \Hom(Y,\CExp(E,X)) &\cong \Hom(Y \times E, X) \cong \Hom(\sk_n(Y) \times E, X) \\
    &\cong \Hom(\sk_n Y, \CExp(E,X)),
\end{align*}
where we use that the comparison map $\sk_n(Y) \times E \hookrightarrow Y \times E$ is an isomorphism on $n$-skeleta. We conclude that $\CExp(E,X)$ is $n$-coskeletal and hence lean.
\end{proof}

The analogous statement for the internal hom $X^M$ where $M$ is a simplicial set can be reduced to the case above. We will use that there exists a functor $\mathcal{E} \colon \s\Set \to \dd\Set$ such that there are natural isomorphisms $Y \otimes M \cong Y \times \mathcal{E}(M)$ for any simplicial set $M$ and dendroidal set $Y$ (cf. \cite[\S 4.2]{HeutsMoerdijk2020Trees}). Note that the functor $\mathcal{E}$ by construction lands in the full subcategory of closed dendroidal sets $\cd\Set$, hence we also obtain natural isomorphisms $Y \otimes M \cong Y \times \mathcal{E}(M)$ for closed dendroidal sets $Y$. In $\od\Set$ we have natural isomorphisms $Y \otimes M \cong Y \times o_! \mathcal{E}(M)$.

\begin{lemma}\label{lemma:CalEDegreewiseFinite}
Let $M$ be a degreewise finite simplicial set. Then $\mathcal{E}(M)$ is a degreewise finite dendroidal set. Moreover, if $M \to N$ is a map of simplicial sets that is an isomorphism on $n$-skeleta, then $\mathcal{E}(M) \to \mathcal{E}(N)$ is an isomorphism on $n$-skeleta.
\end{lemma}

\begin{proof}
The functor $\mathcal{E} \colon \s\Set \to \dd\Set$ is fully determined by its definition on standard simplices
\[\mathcal{E}(\Delta[n])_T := \Hom(E(T),[n]) \]
and the fact that it preserves colimits. Here $E(T)$ denotes the poset of edges of the tree $T$, where the partial order is defined by $e \leq f$ if and only if the unique path from $e$ to the root of the tree $T$ passes through $f$. By $\Hom(E(T),[n])$, we mean the set of order-preserving maps $E(T) \to [n]$.

Clearly, $\mathcal{E}(\Delta[n])$ is degreewise finite for every $n$. If $M$ is both degreewise finite and skeletal, then it follows that $\mathcal{E}[M]$ is degreewise finite by writing $M$ as a finite colimit of representables. For the general case, note that $M = \colim \sk_n M$ and that for every $n$, the map $\sk_{n-1} M \to \sk_n M$ is a pushout of a map of the form $\coprod_i \partial \Delta[n] \to \coprod_i \Delta[n]$. Now observe that $\mathcal{E}(\partial \Delta[n])$ can be identified with the subobject of $\mathcal{E}(\Delta[n])$ consisting precisely of the non-surjective maps $E(T) \to [n]$. If $|T| < n$, then $\# E(T) < n+1$, hence $\mathcal{E}(\partial \Delta[n])_T = \mathcal{E}(\Delta[n])_T$. Since $\mathcal{E}$ preserves coproducts and pushouts, this implies that $\mathcal{E}(\sk_{n-1} M)_T = \mathcal{E}(\sk_n M)_T$ whenever $|T| < n$. In particular, we see that $\mathcal{E}(M)_T = \mathcal{E}(\sk_{|T|} M)_T$. It follows that that $\mathcal{E}(M)$ is degreewise finite since this holds for $\mathcal{E}(\sk_n M)$ for every $n \in \bbN$.

The second statement of the lemma follows by the same argument.
\end{proof}

\begin{lemma}\label{lemma:TrueExponentialLean}
Let $Y$ be a lean dendroidal set. Then for any degreewise finite simplicial set $M$, the cotensor $Y^M := Y^{i_! M}$ is again lean, and for any degreewise finite dendroidal set $X$, the simplicial hom $\sHom(X,Y)$ is lean.
\end{lemma}

\begin{proof}
By adjunction, it follows that there are natural isomorphisms $Y^M \cong \CExp(\mathcal{E}(M),Y)$. In particular, it follows from Lemmas \ref{lemma:CartesianExponentialLean} and \ref{lemma:CalEDegreewiseFinite} that $Y^M$ is lean for any degreewise finite simplicial set $M$ and any lean dendroidal set $Y$.

To see that $\sHom(X,Y)$ is lean whenever $X$ is degreewise finite and $Y$ is lean, let $m$ be such that $Y$ is $m$-coskeletal. Note that
\begin{align*}
    \sHom(X,Y)_n &\cong \Hom(X \times \mathcal{E}(\Delta[n]), Y) \cong \Hom(X \times \mathcal{E}(\sk_m \Delta[n]), Y) \\
    &\cong \Hom(\sk_m \Delta[n], \sHom(X,Y))
\end{align*}
by adjunction, where we use that $\mathcal{E}(\Delta[n])$ and $\mathcal{E}(\sk_m \Delta[n])$ agree on their $m$-skeleta by \autoref{lemma:CalEDegreewiseFinite}. From this description it follows that $\sHom(X,Y)$ is $m$-coskeletal and degreewise finite, hence lean. The cases of closed and open dendroidal sets are similar, using $o_!\mathcal{E}$ instead of $\mathcal{E}$ in the latter case.
\end{proof}

\subsection{Lean \texorpdfstring{$\infty$}{infinity}-operads} \label{ssec:LeanInfinityOperads}

The goal of this section is to give a homotopy-theoretic characterization of which $\infty$-operads are equivalent to lean $\infty$-operads in \autoref{theorem:PiFiniteInftyOperadVsLean} below. To motivative our results, we first sketch the proof that a simplicial set is weakly equivalent to a lean Kan complex if and only if it is \emph{$\pi$-finite}. The latter is also proved in \cite[Lemma 7.2.4]{BarneaHarpazHorel2017} in a slightly different way. (Note that what we call lean is called \emph{$\tau$-finite} there.)

\paragraph{Lean Kan complexes.}

\begin{definition}\label{def:PiFiniteSimplicialSets}
	A simplicial set $X$ is called \emph{$\pi$-finite} if $\pi_n(X,x_0)$ is finite for each $n \geq 0$ and each vertex $x_0$, and if moreover there exists an $n_0$ such that $\pi_n(X,x_0) = \{0\}$ for each $n \geq n_0$ and each vertex $x_0$.
\end{definition}

\begin{proposition}\label{prop:BasicPropertiesPiFiniteKan}
\begin{enumerate}[(i)]
    \item\label{item1:prop:BasicPropertiesPiFiniteKan} If $X \wearrow Y$ is a weak equivalence, then $X$ is $\pi$-finite if and only if $Y$ is.
    \item\label{item2:prop:BasicPropertiesPiFiniteKan} If $p \colon E \fibarrow X$ is a fibration with $X$ $\pi$-finite, then $E$ is $\pi$-finite if and only if the fiber $E|_x = p^{-1}(x)$ is $\pi$-finite for any vertex $x$ of $X$.
    \item\label{item3:prop:BasicPropertiesPiFiniteKan} Any finite homotopy limit of $\pi$-finite simplicial sets is again $\pi$-finite.
\end{enumerate}
\end{proposition}

\begin{proof}
Item \ref{item1:prop:BasicPropertiesPiFiniteKan} is trivial, and item \ref{item2:prop:BasicPropertiesPiFiniteKan} follows by analysing the long exact sequence of a fibration. For \ref{item3:prop:BasicPropertiesPiFiniteKan}, it suffices to show that the the terminal object is $\pi$-finite, and that the class of $\pi$-finite spaces is closed under pullbacks along fibrations and cotensors by the standard simplices $\Delta[n]$. The latter follows from items \ref{item2:prop:BasicPropertiesPiFiniteKan} and \ref{item1:prop:BasicPropertiesPiFiniteKan}, respectively.
\end{proof}

Recall that a Kan complex $X$ is called \emph{minimal} if for any $n \in \bbN$ and for any pair of $n$-simplices $x,y \colon \Delta[n] \to X$ that are homotopic relative to $\partial \Delta[n]$, one has that $x = y$. Any Kan complex contains a minimal Kan complex as a deformation retract (see e.g. \cite[\S VI.5.2]{GabrielZisman1967Calculus}).

One can rephrase the definition of minimality as follows. Given a map $D \colon \partial \Delta[n] \to X$, let $\Fill(D)$ denote the fiber of $X^{\Delta[n]} \to X^{\partial \Delta[n]}$ above $D$. Then $X$ is minimal if and only if for any $n \in \bbN$ and any $D \colon \partial \Delta[n] \to X$, the map $\Fill(D)_0 \to \pi_0(\Fill(D))$ is an isomorphism.

\begin{lemma}\label{lemma:MinimalKanDegreewiseFinite}
Suppose that $X$ is a minimal Kan complex that is $\pi$-finite. Then $X$ is degreewise finite.
\end{lemma}

\begin{proof}
Note that if $X$ has finitely many simplices in degree $< n$, then there are finitely many maps $D \colon \partial \Delta[n] \to X$. In particular, the result follows by induction if we can show that for any boundary $D \colon \partial \Delta[n] \to X$, there are finitely many fillers $\Delta[n] \to X$; that is, if $\Fill(D)$ has finitely many $0$-simplices. By minimality, $\Fill(D)_0 = \pi_0(\Fill(D))$. It follows from \autoref{prop:BasicPropertiesPiFiniteKan} that $X^{\Delta[n]}$ and $X^{\partial \Delta[n]}$ are $\pi$-finite, hence that the fiber $\Fill(D)$ is $\pi$-finite. In particular, $\pi_0(\Fill(D)) = \Fill(D)_0$ is finite.
\end{proof}

\begin{lemma}\label{lemma:TruncatedKanEquivalentToCoskeleton}
Let $X$ be a Kan complex and suppose that $n_0$ is given such that $\pi_n(X,x_0) = \{0\}$ for all vertices $x_0$ and all $n > n_0$. Then $X \to \cosk_{n_0 + 1} X$ is a weak equivalence.
\end{lemma}

\begin{proof}
It is clear that $\cosk_{n_0+1} X$ is a Kan complex and that $X \to \cosk_{n_0+1} X$ induces isomorphisms on homotopy groups of dimension $\leq n_0$. Since $\Hom(\Delta[n],\cosk_{n_0+1} X) \cong \Hom(\partial \Delta[n], \cosk_{n_0+1} X)$ for $n > n_0+1$, the homotopy groups of $\cosk_{n_0+1} X$ above dimension $n_0$ vanish.
\end{proof}

\begin{proposition}\label{prop:PiFiniteKanVsLean}
A simplicial set $X$ is $\pi$-finite if and only if it is weakly equivalent to a lean Kan complex.
\end{proposition}

\begin{proof}
It is clear that any lean Kan complex is $\pi$-finite. For the converse, let $M$ be any minimal Kan complex weakly equivalent to $X$. By Lemmas \ref{lemma:MinimalKanDegreewiseFinite} and \ref{lemma:TruncatedKanEquivalentToCoskeleton}, $M$ is degreewise finite and there exists an $n$ such that $M \simeq \cosk_{n} M$. Since $\cosk_{n} M$ is a lean Kan complex, the result follows.
\end{proof}

\begin{remark}
In the above proof, it is not necessary to replace the minimal Kan complex $M$ by $\cosk_n M$, since any minimal Kan complex whose homotopy groups vanish above a certain dimension can be shown to be coskeletal.
However, the analogous statement for dendroidal sets turns out to be false: if $X$ is a minimal $\infty$-operad equivalent to a lean $\infty$-operad, then $X$ need not be coskeletal.
\end{remark}

\paragraph{Lean $\infty$-operads.} 

Recall that an (open) $\infty$-operad is defined as a fibrant object in the operadic model structure on the category of (open) dendroidal sets $X$. Throughout what follows, all results hold for both open and general dendroidal sets unless stated otherwise, so we will drop the adjective ``open''. The case of closed dendroidal sets is somewhat different and will be discussed below. Throughout the rest of this section, to avoid cluttered notation we will write $T$ for the representable presheaf $\Omega[T]$ that a tree $T$ represents. Recall from \autoref{ssec:dSetsAndSpaces} the definitions of the mapping space $\Map(X,Y)$ and the spaces of operations of an $\infty$-operad.

\begin{definition}\label{def:PiFiniteInftyOperad}
An $\infty$-operad $X$ is called \emph{$\pi$-finite} if
\begin{enumerate}[(i)]
    \item\label{item1:def:PiFiniteInftyOperad} the set $\pi_0(k i^* X) = \pi_0(\Map(\eta,X))$ is finite,
    \item\label{item2:def:PiFiniteInftyOperad} for any tuple of colours $c_1,\ldots,c_n,d \in X_\eta$, the space of operations $X(c_1,\ldots,c_n;d)$ is $\pi$-finite in the sense of \autoref{def:PiFiniteSimplicialSets}, and
    \item\label{item3:def:PiFiniteInftyOperad} there exists an $n_0$ such that for any $n > n_0$ and any tuple of colours $c_1,\ldots,c_n,d \in X_\eta$, the space of operations $X(c_1,\ldots,c_n;d)$ is a contractible Kan complex.
\end{enumerate}
\end{definition}

\begin{remark}\label{remark:PiFiniteInftyOperadTruncatedHomotopy}
Item \ref{item1:def:PiFiniteInftyOperad} can be rephrased by saying that $X$ has finitely many colours up to equivalence. Namely, two colours $c$ and $d$ represent the same element in $\pi_0(k i^* X)$ if and only if there are maps $c \to d$ and $d \to c$ in the underlying $\infty$-category $i^*X$ of $X$ that are inverse to each other.
\end{remark}

\begin{remark}
By combining all three items of this definition, it follows that for any $\pi$-finite $\infty$-operad $X$, there exists an $m_0$ such that the homotopy groups of all the spaces of operations of $X$ vanish above dimension $m_0$.
\end{remark}

Since weak equivalences of $\infty$-operads are essentially surjective and fully faithful, it follows that if $X \wearrow Y$ is a weak equivalence of $\infty$-operads, then $X$ is $\pi$-finite if and only if $Y$ is.

The following is the main result of this section.

\begin{theorem}\label{theorem:PiFiniteInftyOperadVsLean}
An $\infty$-operad is $\pi$-finite if and only if it is weakly equivalent to a lean $\infty$-operad.
\end{theorem}

The proof will require a few results on minimal $\infty$-operads and coskeleta of $\infty$-operads. The definition of a minimal Kan complex is generalized to $\infty$-operads as follows. We say that  $x,y \colon T \to X$ are \emph{homotopic relative to their boundary} if they agree on $\partial T$ and there exists a lift in
\[\begin{tikzcd}
\{0,1\} \ar[d,hook] \ar[rr,"(x{,}y)"] & &[50pt] \sHom(T,X) \ar[d] \\
J \ar[r] \ar[urr, dashed] & \{*\} \ar[r, "x|_{\partial T} = y|_{\partial T}"'] & \sHom(\partial T,X),
\end{tikzcd}\]
where $J$ is the nerve of the groupoid $(0 \leftrightarrow 1)$. The following reformulation will be convenient for us. Write $\Fill(D)$ for the fiber of $\Map(T,X) \fibarrow \Map(\partial T,X)$ above a boundary $D \colon \partial T \to X$.

\begin{lemma} \label{lem:HoRelBoundaryAndFill}
Let $X$ be an $\infty$-operad and $x,y \in X_T$ two elements that agree on their boundary $D$. Then $x$ and $y$ are homotopic relative to their boundary if and only if they are in the same path component of $\Fill(D)$.
\end{lemma}

\begin{proof}
Since $J$ is a Kan complex, any map from $J$ to $\sHom(T,X)$ lands in $k\sHom(T,X) = \Map(T,X)$. In particular, we can replace $\sHom$ by $\Map$ in the above definition of ``homotopic relative to their boundary''. This implies that two elements $x,y \in X_T$ with $x|_{\partial T} = y|_{\partial T} = D$ are homotopic relative to their boundary if and only if there exists a map $f \colon J \to \Fill(D)$ such that $f(0) = x$ and $f(1) = y$. Since $\Delta[1] \cofarrow J$ is a trivial cofibration in the Kan-Quillen model structure on $\s\Set$, this is the case if and only if $x$ and $y$ are in the same path component of $\Fill(D)$.
\end{proof}

An $\infty$-operad is called \emph{minimal} if for any $x,y \in X_T$ that are homotopic relative to their boundary, there exists an $\alpha \in \Aut(T)$ such that $y = \alpha^*(x)$. It can be shown that any normal $\infty$-operad $X$ contains a minimal $\infty$-operad $M \trivcofarrow X$ as a deformation retract (see \cite[Theorem 1.1]{MoerdijkNuiten2016Minimal}).

\begin{proposition}\label{prop:PiFiniteMappingSpacesForOperads}
Let $X$ be a $\pi$-finite $\infty$-operad. Then for any tree $T$, the Kan complex $\Map(T,X)$ is $\pi$-finite, and for any boundary $D \colon \partial T \to X$, the Kan complex $\Fill(D)$ is $\pi$-finite.
\end{proposition}

\begin{proof} First consider the case $T = \eta$. The set $\pi_0\Map(\eta,X)$ is finite by definition. Since $\Map(\eta,X)$ is the maximal Kan complex contained in the underlying $\infty$-category $i^*X$ of $X$, the loop space of $\Map(\eta,X)$ at $c \in X_\eta$ agrees with the subspace of $X(c;c)$ that is the union of all path components whose $0$-vertices are equivalences from $c$ to itself in $i^*X$. By item \ref{item2:def:PiFiniteInftyOperad} of \autoref{def:PiFiniteInftyOperad}, this space is $\pi$-finite.

Now consider the case $T = C_n$. Since $\sqcup_{0 \leq i \leq n} \eta = \partial C_n$, we see that $\Map(\partial C_n, X)$ is a $\pi$-finite Kan complex. The fibration $\Map(C_n,X) \fibarrow \Map(\partial C_n,X)$ has $\pi$-finite fibers $X(c_1,\ldots,c_n;d)$ by item \ref{item2:def:PiFiniteInftyOperad} of \autoref{def:PiFiniteInftyOperad}, hence $\Map(C_n,X)$ is $\pi$-finite by \autoref{prop:BasicPropertiesPiFiniteKan}.

Now suppose that $T$ is a tree with at least two vertices. By Lemma 6.37 of \cite{HeutsMoerdijk2020Trees}, the inclusion $\Spine[T] \cofarrow T$ is inner anodyne, where $\Spine[T]$ denotes the union of all external faces of $T$ with exactly one vertex, hence $\Map(T,X) \trivfibarrow \Map(\Spine[T],X)$ is a weak equivalence. $\Spine[T]$ is a finite (homotopy) colimit of faces of the form $C_n$ and $\eta$, hence $\Map(\Spine[T],X)$ is a finite homotopy limit of Kan complexes of the form $\Map(C_n,X)$ and $\Map(\eta,X)$. By \autoref{prop:BasicPropertiesPiFiniteKan}, $\Map(T,X)$ is $\pi$-finite.

Finally, suppose that $D \colon \partial T \to X$ is given for some tree $T$. Note that $\partial T$ is the finite (homotopy) colimit of all proper faces of $T$, hence $\Map(\partial T, X)$ is the finite homotopy limit of the Kan complexes $\Map(S,X)$ where $S$ ranges over all proper faces of $T$. By \autoref{prop:BasicPropertiesPiFiniteKan}, both $\Map(\partial T, X)$ and the fiber $\Fill(D)$ of $\Map(T,X) \fibarrow \Map(\partial T, X)$ must be $\pi$-finite.
\end{proof}

\begin{corollary}\label{cor:PiFiniteMinimalInftyOperadDegwiseFinite}
Suppose that $X$ is a minimal $\pi$-finite $\infty$-operad. Then $X$ is degreewise finite.
\end{corollary}

\begin{proof}
If $X$ is minimal, then by \autoref{lem:HoRelBoundaryAndFill} we see that for any boundary $D \colon \partial T \to X$, the number of fillers $T \to X$ is at most $|\pi_0(\Fill(D)) \times \Aut(T)| < \infty$. The result now follows from \autoref{prop:PiFiniteMappingSpacesForOperads} by induction.
\end{proof}

To study the behaviour of the coskeleton functors with respect to $\infty$-operads, we need the following two lemmas.

\begin{lemma}\label{lemma:SpacesOfOperationsOfHomotopyCoherentNerve}
	Let $\mathcal{P}$ be a fibrant simplicial operad. Then the spaces of operations of the homotopy coherent nerve $w^* \mathcal{P}$ are canonically equivalent to those of $\mathcal{P}$.
\end{lemma}

\begin{proof}
	Let $\wh{\mathcal{P}}_\bullet$ be a Reedy fibrant simplicial resolution of $\mathcal{P}$, where we may assume that $\wh{\mathcal{P}}_0 = \mathcal{P}$ since $\mathcal{P}$ is fibrant, and let $(c_1,\ldots,c_n,d)$ be a sequence of colours in $\mathcal P$. We can view this as a sequence of colours in $\wh{\mathcal P}_m$ for every $m$ through the degeneracy maps of $\wh{\mathcal P}_\bullet$. Recall that $(w^* \mathcal P)(c_1,\ldots,c_n;d)$ is the fiber of $\Map(C_n, w^* \mathcal P) \fibarrow \Map(\eta, w^* \mathcal P)^{n+1}$, which is equivalent to the fiber of $\Hom(w_! C_n, \wh{\mathcal P}_\bullet) \to \Hom(w_! \eta, \wh{\mathcal P}_\bullet)^{n+1}$ above $(c_1,\ldots,c_n, d)$ via the Quillen equivalence $w_! \dashv w^*$. This fiber is easily shown to agree with the Kan complex $\Hom(\Delta[0], \wh{\mathcal P}_\bullet(c_1,\ldots,c_n;d))$. Since $\wh{\mathcal P}_\bullet(c_1,\ldots,c_n;d)$ is a simplicial resolution of $\mathcal{P}(c_1,\ldots,c_n;d)$, this Kan complex must be homotopy equivalent to $\mathcal P(c_1,\ldots,c_n;d)$.
\end{proof}

\begin{lemma}\label{lem:RlpWrtBoundaryInclusionsLargeTrees}
Let $X$ be a $\pi$-finite $\infty$-operad. Then there exists a $q_0$ such that for all trees $T$ of size $|T| > q_0$, the $\infty$-operad $X$ has the \rlp \wrt $\partial T \cofarrow T$.
\end{lemma}

\begin{proof}
Note that if $X \trivfibarrow Y$ is a trivial fibration, then since $\partial T$ is normal (i.e. cofibrant), a standard argument shows that $X$ has the \llp \wrt $\partial T \to T$ if and only if $Y$ does. By the Quillen equivalence
\[w_! : \dd\Set \rightleftarrows \s\Op : w^*,\]
any $\infty$-operad $X$ is equivalent to one of the form $w^* \mathcal{P}$, where $\mathcal{P}$ is a fibrant simplicial operad. An application of Brown's lemma now shows that $X$ has the \llp \wrt $\partial T \to T$ if and only if $w^* \mathcal{P}$ does (see e.g. \cite[Corollary 7.7.2]{Hirschhorn2003Model}). Throughout the rest of the proof, assume without loss of generality that $X = w^* \mathcal{P}$.

Since the $\infty$-operad $X = w^*\mathcal{P}$ is $\pi$-finite, by \autoref{lemma:SpacesOfOperationsOfHomotopyCoherentNerve} there exist $m_0, n_0 \in \bbN$ such that for any tuple of colours $c_1,\ldots,c_n,d$ of $\mathcal P$ and any vertex $p \in \mathcal{P}(c_1,\ldots,c_n;d)$, one has $\pi_m(\mathcal{P}(c_1,\ldots,c_n;d),p) = 0$ whenever $m > m_0$ or $n > n_0$. A simple counting argument shows that there exists a $q_0$ such that any tree $T$ with $|T| > q_0$ has more than $m_0 + 1$ internal edges or more than $n_0$ leaves. We will show that $w^*\mathcal{P}$ has the \rlp \wrt $\partial T \cofarrow T$ for any such tree $T$.

By adjunction, we need to construct a lift in
\[\begin{tikzcd}
w_! \partial T \ar[d] \ar[r,"f"] & \mathcal{P} \\
w_! T \ar[ur, dashed] &
\end{tikzcd}\]
An inspection of the functor $w_! \colon \dd\Set \to \s\Op$ as in the proof of \cite[Proposition 4.5]{CisinskiMoerdijk2013SimplicialOperads} shows that $(w_! \partial T)(c_1,\ldots,c_n;d) \to (w_! T)(c_1,\ldots,c_n;d)$ is always an isomorphism of simplicial sets, unless $c_1,\ldots,c_n$ is exactly the set of leaves of $T$ and $d$ the root of $T$, in which case it is the inclusion of the boundary of a cube $\partial(\Delta[1]^{\mathrm{in}(T)}) \hookrightarrow \Delta[1]^{\mathrm{in}(T)}$ where $\mathrm{in}(T)$ denotes the number of internal edges of $T$. In particular, constructing a lift in the diagram above comes down to constructing a lift in
\[\begin{tikzcd}
\partial(\Delta[1]^{\mathrm{in}(T)}) \ar[r] \ar[d] & \mathcal{P}(f(e_1),\ldots,f(e_n);f(r)) \\
\Delta[1]^{\mathrm{in}(T)} \ar[ur, dashed] &
\end{tikzcd}\]
where $e_1,\ldots,e_n$ are the leaves of $T$ (in an arbitrary but fixed order) and $r$ is the root edge of $T$. The obstruction to such a lift lies in $\pi_{\mathrm{in}(T)-1}(\mathcal{P}(f(e_1),\ldots,f(e_n);f(r)))$, which is trivial by our assumptions on the size of $T$. We conclude that $X$ has the \rlp \wrt $\partial T \cofarrow T$.
\end{proof}

\begin{corollary}\label{lem:PiFiniteInftyOperadEquivToCoskeleton}
Let $X$ be a $\pi$-finite $\infty$-operad. Then for sufficiently large $q$, the unit $X \to \cosk_q X$ is a trivial fibration and $\cosk_q X$ is an $\infty$-operad.
\end{corollary}

\begin{proof}
By adjunction, $X \to \cosk_q X$ has the \rlp \wrt $\partial T \cofarrow T$ if and only if $X$ has the \rlp \wrt $\partial T \cup_{\sk_q \partial T} \sk_q T \to T$. If $|T| \leq q$, then $\sk_q T = T$, so this map is an isomorphism. If $|T| > q$, then $\sk_q T = \sk_q \partial T$, hence this map agrees with the boundary inclusion $\partial T \cofarrow T$. In particular, it follows from \autoref{lem:RlpWrtBoundaryInclusionsLargeTrees} that for $q$ sufficiently large, $X \to \cosk_q X$ is a trivial fibration.

To see that $\cosk_q X$ is fibrant, note that since $X \trivfibarrow \cosk_q X$ is a trivial fibration and any horn $\Lambda^e[T]$ is cofibrant, $\cosk_q X$ has the \rlp \wrt $\Lambda^e[T] \cofarrow T$ if and only if $X$ does.
\end{proof}

\begin{remark}
If $X$ is an open $\infty$-operad (not necessarily $\pi$-finite), then the coskeleton $\cosk_q X$ can be shown to be an (open) $\infty$-operad for any $q \geq 0$. However, this is not true for general $\infty$-operads.
\end{remark}

We are now ready to prove the main result.

\begin{proof}[Proof of \autoref{theorem:PiFiniteInftyOperadVsLean}]
We first show that any lean $\infty$-operad is $\pi$-finite. Let $X$ be a degreewise finite and $q$-coskeletal $\infty$-operad. Clearly $X$ has finitely many colours. Furthermore, by \autoref{lemma:TrueExponentialLean}, we see that $\sHom(\partial C_n, X)$ and $\sHom(C_n, X)$ are lean for every $n \geq 0$, hence the spaces of operations $X(c_1,\ldots,c_n;d)$ are lean. By \autoref{prop:PiFiniteKanVsLean}, we see that $X(c_1,\ldots,c_n;d)$ is a $\pi$-finite Kan complex for every $n$. Finally, for $n \geq q$ we have $|C_n| = q+1$, hence $\sHom(C_n,X) \to \sHom(\partial C_n, X)$ is an isomorphism. In particular, $X(c_1,\ldots,c_n;d) = *$ for $n \geq q$, so we conclude that $X$ is a $\pi$-finite $\infty$-operad.

For the converse, suppose that $X$ is any $\pi$-finite $\infty$-operad. Choose a normalization $\wt X \trivfibarrow X$ and a minimal $\infty$-operad $M \trivcofarrow X$. Then $M$ is degreewise finite by \autoref{cor:PiFiniteMinimalInftyOperadDegwiseFinite}, so by \autoref{lem:PiFiniteInftyOperadEquivToCoskeleton} we see that $\cosk_q M$ is a lean $\infty$-operad weakly equivalent to $X$ for $q$ sufficiently large.
\end{proof}

A similar (but simpler) proof gives the following analogue for $\infty$-categories (cf. \cite[Lemmas A.9-10]{BlomMoerdijk2020SimplicialProV1}).

\begin{proposition}
An $\infty$-category $X$ is equivalent to a lean $\infty$-category if and only if
\begin{enumerate}[(i)]
    \item the set $\pi_0(kX)$ is finite, and
    \item for any two objects $x,y \in X_0$, the mapping space $\map_X(x,y)$ is $\pi$-finite in the sense of \autoref{def:PiFiniteSimplicialSets}.
\end{enumerate}
\end{proposition}

\paragraph{Lean closed $\infty$-operads.}

For a closed tree $T$, denote the representable closed dendroidal set $\Omegacl[T]$ by $T$ and the closed boundary $\clboun \Omegacl[T]$ by $\clboun T$. For open and general dendroidal sets, the boundary $\partial C_n$ is the disjoint union $\sqcup_{0 \leq i \leq n} \eta$ of the root and leaf edges of $C_n$. In the closed case, for $n \geq 1$ one similarly has inclusions $\sqcup_{0 \leq i \leq n} \ol \eta \hookrightarrow \clboun \ol C_n \hookrightarrow \ol C_n$, but the first of these is generally not the identity.

For a closed $\infty$-operad $X$, these inclusions give the diagram of categorical fibrations
\[\begin{tikzcd}[column sep = tiny]
\sHom(\ol C_n, X) \ar[rr, two heads] \ar[dr,two heads] & & \sHom(\clboun \ol C_n, X) \ar[dl, two heads] \\
 & \sHom(\ol \eta, X)^{n+1} &
\end{tikzcd}\]
Recall from \autoref{ssec:dSetsAndSpaces} that the fibers of the left-hand map are the spaces of operations $X(c_1,\ldots,c_n;d)$. The fiber above $(c_1,\ldots,c_n;d)$ of the right-hand map will be denoted $X^-(c_1,\ldots,c_n;d)$ and called the \emph{matching object} of $X$ at $(c_1,\ldots,c_n;d)$. The induced map of fibers $X(c_1,\ldots,c_n;d) \to X^-(c_1,\ldots,c_n;d)$ is the \emph{matching map}. One can show that the inclusion $X^-(c_1,\ldots,c_n;d) \hookrightarrow \sHom(\clboun \ol C_n, X)$ lands in $\Map(\clboun \ol C_n, X) = k \sHom(\clboun \ol C_n, X)$, hence the matching object can also be defined as the fiber of $\Map(\clboun \ol C_n, X) \fibarrow \Map(\ol \eta, X)^n$. In particular, the matching objects $X^-(c_1,\ldots,c_n;d)$ are Kan complexes and the matching maps are Kan fibrations.

If $X$ is lean, then there exists an $n_0$ such that for any $n > n_0$, the map $\sHom(\ol C_n, X) \fibarrow \sHom(\clboun \ol C_n, X)$ is an isomorphism. This implies that for $n > n_0$ and any tuple of colours $c_1,\ldots,c_n,d \in X_{\ol \eta}$, the matching map $X(c_1,\ldots,c_n;d) \fibarrow X^-(c_1,\ldots,c_n;d)$ is an isomorphism. This leads to the following alternative definition of $\pi$-finite in the case of closed $\infty$-operads.

\begin{definition}\label{def:PiFiniteClosedInftyOperad}
A closed $\infty$-operad $X$ is called \emph{$\pi$-finite} if
\begin{enumerate}[(i)]
    \item\label{item1:def:PiFiniteClosedInftyOperad} the set $\pi_0(k i^* X) = \pi_0(\Map(\eta,X))$ is finite,
    \item\label{item2:def:PiFiniteClosedInftyOperad} for any tuple of colours $c_1,\ldots,c_n,d \in X_\eta$, the space of operations $X(c_1,\ldots,c_n;d)$ is $\pi$-finite in the sense of \autoref{def:PiFiniteSimplicialSets}, and
    \item\label{item3:def:PiFiniteClosedInftyOperad} there exists an $n_0$ such that for any $n > n_0$ and any tuple of colours $c_1,\ldots,c_n,d \in X_\eta$, the matching map $X(c_1,\ldots,c_n;d) \fibarrow X^-(c_1,\ldots,c_n;d)$ is a weak equivalence.
\end{enumerate}
\end{definition}

\begin{remark}\label{remark:MatchingObjectAsCubicalLimit}
To gain a better understanding of the matching object of a closed $\infty$-operad $X$, we need to take a closer look at $\clboun \ol C_n$. If $n = 1$, then $\clboun \ol C_1 = \sqcup_{i=0,1} \ol \eta$, hence the matching object $X^-(c;d)$ is a point for any two colours $c,d \in X_{\ol \eta}$. For the case $n > 1$, label the non-root edges of $\ol C_n$ by $1,\ldots, n$. To any proper subset $U \subsetneq \ul n = \{1,\ldots, n\}$ we associate the subobject $\ol C_U \subset \Omegacl[\ol C_n]$ spanned by those faces whose non-root edges are contained in $U$. Clearly $\ol C_U \cong \Omegacl[\ol C_{\# U}]$. Ordering the subsets $U \subsetneq \ul n$ by inclusion, we see that $\clboun \ol C_n = \colim_{U \subsetneq \ul n} \ol C_U$. Since for any $U \subsetneq \ul n$ the map $\colim_{V \subsetneq U} C_V \to C_U$ is a normal monomorphism, this is a homotopy colimit. We deduce that $\Map(\clboun \ol C_n, X) = \lim_U \Map(\ol C_U, X)$ is a finite homotopy limit. Now let a tuple of colours $c_1,\ldots,c_n,d \in X_{\ol \eta}$ be given. For $U \subset \ul n$, let $c_U$ denote the tuple $(c_i)_{i \in U}$. The isomorphisms $\ol C_U \cong \Omegacl[\ol C_{\# U}]$ tell us that the fiber $X^-(c_1,\ldots,c_n;d)$ of
\[\lim_U \Map(\ol C_U, X) = \Map(\clboun \ol C_n, X) \fibarrow \Map(\ol \eta, X)^{n+1}\]
is isomorphic to $\lim_{U \subsetneq \ul n} X(c_U;d)$ , providing a more concrete description of the matching object. Note that this description also holds when $n = 1$. This also justifies the name ``matching object'', since the matching objects of a closed operad $\mathcal{P}$ are defined by $\mathcal{P}^-(c_1,\ldots,c_n;d) = \lim_{U \subsetneq \ul n} \mathcal{P}(c_U;d)$. An argument similar to the one used in \autoref{lemma:SpacesOfOperationsOfHomotopyCoherentNerve} proves that if $\mathcal{P}$ is a Reedy fibrant closed simplicial operad, then $\mathcal{P}(c_1,\ldots,c_n;d) \fibarrow \mathcal{P}^-(c_1,\ldots,c_n;d)$ is canonically equivalent to the matching map $(\ol w^* \mathcal{P})(c_1,\ldots,c_n;d) \fibarrow (\ol w^* \mathcal{P})^-(c_1,\ldots,c_n;d)$ of the homotopy coherent nerve of $\mathcal{P}$.
\end{remark}

\begin{example}\label{exx:AssociativeOperadCoskeletal}
	The (closed nerve of the) unital associative operad $\mathcal{A}ss$ can be shown to be $\pi$-finite in the sense of \autoref{def:PiFiniteClosedInftyOperad}. It is clear that items \ref{item1:def:PiFiniteClosedInftyOperad} and \ref{item2:def:PiFiniteClosedInftyOperad} are satisfied. We will show that the matching map $\mathcal{A}ss(n) \to \mathcal{A}ss^-(n)$ is an isomorphism for $n \geq 4$. Note that for every $n$, one can identify the set $\mathcal{A}ss(n)$ with the set of linear orderings of $\ul n = \{1,\ldots,n\}$. By the previous remark, we can identify the matching object $\mathcal{A}ss^-(n)$ with the limit $\lim_{U \subsetneq \ul n} \mathcal{A}ss(U)$, where we write $\mathcal{A}ss(U)$ for the set of linear orders on the subset $U$ and where for any subset inclusion $V \subset U$, the map $\mathcal{A}ss(U) \to \mathcal{A}ss(V)$ is defined by restricting a linear order on $U$ to one on $V$. It is then easy to verify that for $n \geq 4$, the map $\mathcal{A}ss(n) \to \mathcal{A}ss^-(n)$ is indeed an isomorphism. A similar analysis can be used to show that the closed nerve of $\mathcal{A}ss$ is furthermore a lean closed $\infty$-operad.
\end{example}

\begin{remark}
Call a Kan complex $m_0$-truncated if all its homotopy groups above dimension $m_0$ are trivial. As in \autoref{remark:PiFiniteInftyOperadTruncatedHomotopy}, for any $\pi$-finite closed $\infty$-operad $X$, there exists an $m_0$ such that all the spaces of operations of $X$ are $m_0$-truncated. To see this, let $n_0$ be as in item \ref{item3:def:PiFiniteClosedInftyOperad} of \autoref{def:PiFiniteClosedInftyOperad}. By items \ref{item1:def:PiFiniteClosedInftyOperad} and \ref{item2:def:PiFiniteClosedInftyOperad}, there exists an $m_0$ such that for all $n \leq n_0$, all spaces of $n$-ary operations $X(c_1,\ldots,c_n;d)$ are $m_0$-truncated. It follows inductively that this also holds for $n > n_0$. Namely, suppose that $n > n_0$ is given and that the spaces of $n'$-ary operations are $m_0$-truncated for any $n' < n$. It is well-known that $m_0$-truncated spaces are closed under homotopy limits, hence by \autoref{remark:MatchingObjectAsCubicalLimit} the matching objects $X^-(c_1,\ldots,c_n;d)$ are $m_0$-truncated. By item \ref{item3:def:PiFiniteClosedInftyOperad}, the spaces of $n$-ary operations $m_0$-truncated for all $n$.
\end{remark}

It is clear that if $X \wearrow Y$ is a weak equivalence of closed $\infty$-operads, then $X$ is $\pi$-finite if and only if $Y$ is. We will prove the following analogue of \autoref{theorem:PiFiniteInftyOperadVsLean}.

\begin{theorem}\label{theorem:PiFiniteClosedInftyOperadVsLean}
A closed $\infty$-operad is $\pi$-finite if and only if it is weakly equivalent to a lean closed $\infty$-operad.
\end{theorem}

The proof follows the same steps as that of \autoref{theorem:PiFiniteInftyOperadVsLean}. For any map $D \colon \clboun T \to X$, define the Kan complex $\Fill(D)$ as the fiber of $\Map(T,X) \fibarrow \Map(\clboun T, X)$ above $D$.

\begin{proposition}
Let $X$ be a $\pi$-finite closed $\infty$-operad. Then for any tree $T$, the Kan complex $\Map(T,X)$ is $\pi$-finite, and for any boundary $D \colon \clboun T \to X$, the Kan complex $\Fill(D)$ is $\pi$-finite.
\end{proposition}

\begin{proof}
The proof is similar to that of \autoref{prop:PiFiniteMappingSpacesForOperads}, except that one uses the fibrations $\Map(\ol C_n, X) \fibarrow \Map(\ol \eta, X)^{n+1}$ in place of $\Map(C_n,X) \fibarrow \Map(\partial C_n,X)$ and the \emph{closed} spine of \cite[\S 6]{Moerdijk2018Closed} in place of the usual spine to prove that $\Map(T,X)$ is $\pi$-finite for general trees $T$.
\end{proof}

The theory of minimal $\infty$-operads generalizes to closed $\infty$-operads by \cite[Theorem 5.3]{MoerdijkNuiten2016Minimal}. It follows as in the proof of \autoref{lem:HoRelBoundaryAndFill} that two elements $x,y \in X_T$ are homotopic relative to their boundary if and only if they are in the same path component of $\Fill(D)$. The proof of the following result is identical to that of \autoref{cor:PiFiniteMinimalInftyOperadDegwiseFinite}.

\begin{proposition}
Suppose that $X$ is a minimal $\pi$-finite closed $\infty$-operad. Then $X$ is degreewise finite.
\end{proposition}

The analogue of \autoref{lem:RlpWrtBoundaryInclusionsLargeTrees} can be proved for closed dendroidal sets as well.

\begin{lemma}
Let $X$ be a $\pi$-finite closed $\infty$-operad. Then there exists a $q_0$ such that for all closed trees of size $|T| > q_0$, the closed $\infty$-operad $X$ has the \rlp \wrt $\clboun T \cofarrow T$.
\end{lemma}

\begin{proof}
The proof is similar to that of \autoref{lem:RlpWrtBoundaryInclusionsLargeTrees}, with the differences that one instead uses the Quillen equivalence
\[\ol w_! : \cd\Set \rightleftarrows \mathbf{c}\s\Op : \ol w^*\]
to replace $X$ by $\ol w^* \mathcal{P}$ for some Reedy fibrant closed simplicial operad $\mathcal{P}$, and that solving a lifting problem of the form
\[\begin{tikzcd}
\ol w_!(\clboun T) \ar[r,"f"] \ar[d] & \mathcal{P} \\
\ol w_!(T) \ar[ur, dashed] &
\end{tikzcd}\]
is equivalent to solving one of the form
\begin{equation}\label{diag:Closedw!Lifting}
\begin{tikzcd}
\partial(\Delta[1]^{\mathrm{in}(T)}) \ar[r] \ar[d,tail] & \mathcal{P}(f(e_1),\ldots,f(e_n);f(r)) \ar[d, two heads] \\
\Delta[1]^{\mathrm{in}(T)} \ar[r] \ar[ur,dashed] & \mathcal{P}^-(f(e_1),\ldots,f(e_n);f(r)),
\end{tikzcd}
\end{equation}
as shown in the proof of \cite[Theorem 10.4]{Moerdijk2018Closed}. Here $e_1, \ldots, e_n$ are the edges of $T$ immediately below the stumps (in any fixed order), $r$ is the root of $T$ and $\mathrm{in}(T)$ denotes the number very inner edges (the edges not equal to one of $e_1, \ldots, e_n,r$). The obstruction to finding such a lift lies in $\pi_{\mathrm{in}(T)-1}(F)$, where $F$ is the (homotopy) fiber of the right-hand map of diagram \eqref{diag:Closedw!Lifting}. By \autoref{remark:MatchingObjectAsCubicalLimit}, this map is equivalent to the matching map $(\ol w^*\mathcal{P})(f(e_1),\ldots,f(e_n);f(r)) \fibarrow (\ol w^*\mathcal{P})^-(f(e_1),\ldots,f(e_n);f(r))$, hence their fibers are homotopy equivalent. Since $\ol w^*\mathcal{P}$ is $\pi$-finite, by the long exact sequence of homotopy groups there exist $m_0,n_0 \in \bbN$ such that $\pi_{\mathrm{in}(T)-1}(F)$ is trivial whenever $\mathrm{in}(T) - 1 > m_0$ or $n > n_0$. In particular, we see that $X$ has the \rlp \wrt $\clboun T \cofarrow T$ when $|T|$ is sufficiently large.
\end{proof}

\begin{corollary}
Let $X$ be a $\pi$-finite closed $\infty$-operad. Then for sufficiently large $q$, the map $X \to \cosk_q X$ is a trivial fibration and $\cosk_q X$ is a closed $\infty$-operad.
\end{corollary}

\begin{proof}
The proof is analogous to that of \autoref{lem:PiFiniteInftyOperadEquivToCoskeleton}.
\end{proof}

\begin{remark}
As in the open case, one can in fact show that for any closed $\infty$-operad $X$ and every $q \geq 0$, the coskeleton $\cosk_q X$ is again a closed $\infty$-operad.
\end{remark}

\begin{proof}[Proof of \autoref{theorem:PiFiniteClosedInftyOperadVsLean}]
Given the results above, the proof is almost identical to that of \autoref{theorem:PiFiniteInftyOperadVsLean} and left to the reader.
\end{proof}

\section{Pro-categories} \label{sec:ProCategories}

In this section we introduce some conventions and notation regarding pro-categories and prove a few elementary properties of profinite sets.

\subsection{Preliminaries on pro-categories} \label{ssec:ProCategories}

To begin with, we present some basic properties of pro-categories that are used throughout the rest of this paper. Most of these will be familiar to the reader, so we will not give proofs. For details, we refer to \cite[Exposé 1]{grothendieck1972theorie}, \cite[\S 2.1]{EdwardsHastings1976CechSteenrod}, \cite[Appendix]{ArtinMazur1986Etale}, \cite[\S VI.1]{Johnstone1982StoneSpaces} and \cite{Isaksen2002Calculating}. In the discussion below, all (co)limits are asssumed to be \textbf{small}.

For a category $\bfC$, its pro-completion $\Pro(\bfC)$ is obtained by freely adjoining cofiltered (or codirected) limits to $\bfC$. By a \emph{cofiltered limit}, we mean the limit of a diagram indexed by a cofiltered category, i.e. a category $I$ such that $I^{op}$ is filtered. Similarly, a \emph{codirected} limit is a limit indexed by a codirected set $I$. One way to make the idea of ``freely adjoining cofiltered limits'' precise, is to define the objects in $\Pro(\bfC)$ to be all diagrams $I \to \bfC$ for all cofiltered categories $I$. Such objects are called \emph{pro-objects} and denoted $C = \{C_i\}_{i \in I}$. The morphisms between two such objects $C = \{C_i\}_{i \in I}$ and $D = \{D_j\}_{j \in J}$ are defined by
\begin{equation}\label{equation:DefinitionHomSetProCategory}
    \Hom_{\Pro(\bfC)}(C,D) = \lim_j \colim_i \Hom_{\bfC}(C_i,D_j).
\end{equation}

There is a fully faithful embedding $\bfC \hookrightarrow \Pro(\bfC)$ sending an object $C$ to the constant diagram $* \to \bfC$ with value $C$. We will generally identify $\bfC$ with its image in $\Pro(\bfC)$ under this embedding, and abusively write $C$ for the image of an object $C$ under the functor $\bfC \hookrightarrow \Pro(\bfC)$. This embedding preserves all colimits and all finite limits that exist in $\bfC$. Moreover, for any pro-object $\{C_i\}_{i \in I}$ given by a diagram $I \to \bfC$, the cofiltered limit of the composition $I \to \bfC \hookrightarrow \Pro(\bfC)$ agrees with the pro-object $\{C_i\}_{i \in I}$; that is, any pro-object is ``its own limit'' in $\Pro(\bfC)$.

If $C = \{C_i\}_{i \in I}$ is a pro-object indexed by a codirected set, then for any $j \leq i$ in $I$, we will denote the bonding map by $\pi_{ij} \colon X_j \to X_i$. For any $i \in I$, we write $\pi_i \colon C \to C_i$ for the projection map.

The category $\Pro(\bfC)$ can be characterized by the following universal property: it admits all cofiltered limits, and for any category $\bfE$ that admits cofiltered limits, the inclusion $\bfC \hookrightarrow \Pro(\bfC)$ induces an equivalence of categories $\Fun'(\Pro(\bfC),\bfE) \to \Fun(\bfC,\bfE)$, where $\Fun'(\Pro(\bfC),\bfE)$ denotes the category of functors $\Pro(\bfC) \to \bfE$ that preserve cofiltered limits. In particular, any functor $F \colon \bfC \to \bfE$ admits an essentially unique extension $\wt F \colon \Pro(\bfC) \to \bfE$ that preserves cofiltered limits, which is given by $\wt F(\{C_i\}_{i \in I}) \cong \lim_i F(C_i)$.

Since $\{C_i\}_{i \in I}$ is ``its own limit'', it follows that for any cofinal functor $\theta \colon J \to I$, the comparison map $\{C_i\}_{i \in I} \to \{C_{\theta(j)}\}_{j \in J}$ is an isomorphism in $\Pro(\bfC)$. This can be used to reindex pro-objects in more convenient ways. For example, it can be shown that for any cofiltered category $I$, there exists a codirected \emph{poset} $J$ together with a cofinal functor $J \to I$, with the property that for any $j \in J$, the set $J_{> j}$ is finite. In particular, any pro-object is isomorphic to one indexed by such a poset (cf. Proposition 8.1.6 of \cite[Exposé 1]{grothendieck1972theorie}, or Theorem 2.1.6 of \cite{EdwardsHastings1976CechSteenrod} with a correction just after Corollary 3.11 of \cite{BarneaSchlank2015NewModel}). We state this as a lemma for future reference.

\begin{lemma}[\cite{grothendieck1972theorie}]\label{lemma:ReindexPosetSGA}
	Any pro-object is isomorphic to a pro-object indexed by codirected poset $I$ with the property that for every $i \in I$, the set $I_{>i}$ is finite.
\end{lemma}

Other examples of cofinal functors are decreasing monotone functions on codirected sets: any function $\theta \colon I \to I$ satisfying $\theta(i) \leq i$ is cofinal, hence the comparison map $\{C_i\}_{i \in I} \to \{C_{\theta(i)}\}_{i \in I}$ is an isomorphism.\footnote{Note that our use of the term ``cofinal'' is dual to its common usage in order theory: an order-preserving function $I \to J$ between posets is cofinal if and only if $I^{op} \to J^{op}$ is cofinal as a functor.} Such functions are used in the proofs of \autoref{lem:EquivalentCharacterizationsFreeActionOnComplement} and \autoref{lem:NormalMonoHasIncreasinglyNormalRepresentation}.

Recall that if $\bfE$ is a category that has all cofiltered limits, then an object $C$ in $\bfE$ is called \emph{cocompact} if $\Hom_\bfE(-,C)$ sends cofiltered limits to colimits. One can deduce from the definition of the hom-sets \eqref{equation:DefinitionHomSetProCategory} that any object in the image of $\bfC \hookrightarrow \Pro(\bfC)$ is cocompact. There is the following recognition principle for pro-completions, of which we leave the proof to the reader.

\begin{lemma}[Recognition principle]\label{lemma:RecognitionPrincipleInd}
Let $\bfE$ be a category closed under cofiltered limits and let $\bfC \hookrightarrow \bfE$ be a fully faithful functor. If
\begin{enumerate}[(i)]
    \item any object in $\bfC$ is cocompact in $\bfE$, and
    \item any object in $\bfE$ is a cofiltered limit of objects in $\bfC$,
\end{enumerate}
then the canonical extension $\Pro(\bfC) \to \bfE$, coming from the universal property of $\Pro(\bfC)$, is an equivalence of categories.
\end{lemma}

To avoid size issues, we assume from now on that $\bfC$ is an (essentially) small category. The fact that the category $(\Set^{\bfC})^{op}$ is the free completion of $\bfC$ leads to an alternative description of $\Pro(\bfC)$ that is sometimes easier to work with. Namely, we can think of $\Pro(\bfC)$ as the full subcategory of $(\Set^{\bfC})^{op}$ consisting of those objects which are cofiltered limits of representables. If $\bfC$ is small and has finite limits, as will be the case in all of our examples, then these are exactly the functors $\bfC \to \Set$ that send these finite limits of $\bfC$ to limits in $\Set$ (see the dual of Théorème 8.3.3.(v) of \cite[Exposé 1]{grothendieck1972theorie}); i.e.
\[ \Pro(\bfC) \simeq \mathrm{lex}(\bfC,\Set)^{op}, \]
where the right-hand side stands for (the dual of) the category of left exact functors. From this description, one sees immediately that $\Pro(\bfC)$ has all small colimits and that the inclusion $\Pro(\bfC) \to (\Set^{\bfC})^{op}$ preserves these. The category $\Pro(\bfC)$ also has all small limits in this case. Namely, finite products and pullbacks can be computed ``levelwise'' in $\bfC$ as described in \cite[Appendix 4]{ArtinMazur1986Etale}, while cofiltered limits exist as mentioned above. Note however, that while the inclusion $\Pro(\bfC) \to (\Set^{\bfC})^{op}$ preserves cofiltered limits, it does not preserve all limits.

Another consequence of the fact that finite products and pullbacks in $\Pro(\bfC)$ are computed ``levelwise'' is the following: if $\bfE$ is any complete category and $F \colon \bfC \to \bfE$ a functor that preserves finite limits, then its extension $\wt F \colon \Pro(\bfC) \to \bfE$ coming from the universal property mentioned above also preserves finite limits. Since it preserves cofiltered limits by definition, we conclude that it preserves all limits. In fact, more is true. The above description of $\Pro(\bfC)$ as $\mathrm{lex}(\bfC, \Set)^{op}$ allows us to construct a left adjoint $L$ of $\wt F$. Namely, if we define $L(E)(c) := \Hom(E,Fc)$, then $L(E) \colon \bfC \to \Set$ is left exact, hence $L$ defines a functor $\bfE \to \Pro(\bfC)$. Adjointness follows from the Yoneda lemma. We therefore see that, up to unique natural isomorphism, there is a 1-1 correspondence between finite limit preserving functors $\bfC \to \bfE$ and functors $\Pro(\bfC) \to \bfE$ that admit a left adjoint.

An important example of an adjunction obtained in this way is the pro-completion functor. If $\bfE$ is a complete category and $\bfC$ a small full subcategory closed under finite limits, then the inclusion $\bfC \hookrightarrow \bfE$ induces an adjunction
\[\wh{(\cdot)}_{\bfC} :\bfE \rightleftarrows \Pro(\bfC) : U\]
whose left adjoint we call \emph{pro-completion (relative to $\bfC$)} or \emph{pro-$\bfC$ completion}. In many examples, $\bfC$ is the full subcategory of $\bfE$ consisting of objects that are ``finite'' in some sense, and this left adjoint is better known as \emph{profinite completion}. For instance, in the case of groups, this functor $\wh{(\cdot)}_{\Pro} \colon \Grp \to \Pro(\Fin\Grp)$ is the well-known profinite completion functor for groups.

There is also a two-variable version of the above statement: if $\bfC$ and $\bfD$ are small categories that admit finite limits and $\bfE$ is any complete category, then any functor $F \colon \bfC \times \bfD \to \bfE$ that preserves finite limits in both variables extends to a two-variable adjunction. More precisely, it extends to a functor $\wt F \colon \Pro(\bfC) \times \Pro(\bfD) \to \bfE$ such that there exist functors
\[G \colon \Pro(\bfC)^{op} \times \bfE \to \Pro(\bfD) \quad \text{and} \quad H \colon \Pro(\bfD)^{op} \times \bfE \to \Pro(\bfC) \]
together with natural isomorphisms
\[\Hom_{\bfE}(E,\wt F(C,D)) \cong \Hom_{\Pro(\bfD)}(G(C,E),D) \cong \Hom_{\Pro(\bfC)}(H(D,E),C) .\]
In particular, $\wt F$ preserves limits in each variable separately. This construction will be used in \autoref{ssec:WeakEquivalencesDSets} to define a simplicial hom on the category of dendroidal profinite sets.

Let us return to the basic definition \eqref{equation:DefinitionHomSetProCategory} of morphisms in $\Pro(\bfC)$. If $C = \{C_i\}$ and $D = \{D_i\}$ are objects indexed by the same cofiltered category $I$, then any natural transformation with components $f_i \colon C_i \to D_i$ represents a morphism in $\Pro(\bfC)$. Representations of this type will be called \emph{level representations} or \emph{strict representations}. Up to isomorphism, any morphism in $\Pro(\bfC)$ has such a strict representation (see Corollary 3.2 of \cite[Appendix]{ArtinMazur1986Etale}). One can define the notion of a ``level'' diagram or ``strict'' diagram in a similar way. Given an indexing category $K$, a conceptual way of thinking about this is through the canonical functor
\[ L \colon \Pro(\bfC^K) \to \Pro(\bfC)^K. \]
A strict representation of a diagram $D \colon K \to \Pro(\bfC)$ can be thought of as an object of $D'$ of $\Pro(\bfC^K)$ together with a natural isomorphism $D \cong LD'$. In \cite[\S 4]{Meyer1980Approximation}, the following is proved.

\begin{theorem}[\cite{Meyer1980Approximation}]\label{thm:MeyersResult}
Let $K$ be a finite category and $\bfC$ a small category that admits finite limits. Then the canonical functor
\[\Pro(\bfC^K) \to \Pro(\bfC)^K \]
is an equivalence of categories.
\end{theorem}

This shows in particular that, up to isomorphism, any finite diagram in $\Pro(\bfC)$ admits a strict representation if $\bfC$ is small and has finite colimits.

As explained in \cite[Theorem 2.3]{BlomMoerdijk2020SimplicialProV1}, the proof of Proposition 7.4.1 of \cite{BarneaHarpazHorel2017} can be adapted to prove the following extension of Meyer's result.

\begin{theorem}\label{thm:ProValuedPresheavesIsProSkeletalPresheaves}
Let $\bfC$ be a small category that has finite limits, and let $K$ be a small category that can be written as a union of finite full subcategories. Write $\cosk(\bfC^K)$ for the full subcategory of $\bfC^K$ spanned by those functors $K \to \bfC$ that are isomorphic to the right Kan extension of a functor $K' \to \bfC$ for some finite full subcategory $K' \subset K$. Then $\Pro(\cosk(\bfC^K)) \simeq \Pro(\bfC)^K$.
\end{theorem}

Write $\wh \Set = \Pro(\Fin\Set)$ for the category of profinite sets. If we apply the theorem to the categories $\Delta^{op} = \cup_n \Delta^{op}_{\leq n}$ and $\bfC = \Fin\Set$, then we find that $\s\wh\Set \simeq \Pro(\cosk(\s\Fin\Set))$. This is exactly the equivalence of categories proved in Proposition 7.4.1 of \cite{BarneaHarpazHorel2017}.

The main object of study in this paper is the category $\dd\wh\Set$ of dendroidal profinite sets; that is, the category of $\wh\Set$-valued presheaves of the category $\Omega$ of trees. The theorem above gives us an equivalence $\dd\wh\Set \simeq \Pro(\cosk(\dd\Fin\Set))$, where $\cosk(\dd\Fin\Set)$ agrees with the category $\LdSet$ of lean dendroidal sets by \autoref{remark:AlternativeDefinitionLean}. In particular, we can view the category of dendroidal profinite sets both as the pro-category of $\LdSet$ and as the category of $\wh\Set$-valued presheaves on $\Omega$. The notation $\dd\wh\Set$ will refer to either of these categories.

\subsection{Some elementary facts about profinite sets} \label{ssec:ElementaryFactsProfiniteSets}

In this section we will collect some elementary facts about monomorphisms and finite group actions in the category $\wh\Set = \Pro(\Fin\Set)$ of profinite sets, which will be useful when studying normal monomorphisms of dendroidal profinite sets in \autoref{sec:NormalMonos}. It is well-known that $\wh\Set \simeq \Stone$, where $\Stone$ is the category of Stone spaces, i.e. compact Hausdorff totally disconnected spaces. This equivalence in the direction $\wh\Set \to \Stone$ is defined by sending a pro-object $\{X_i\}_{i \in I}$ to its limit $\lim_i X_i$ computed in $\Top$, where each finite set $X_i$ is endowed with the discrete topology. We will interchangeably view profinite sets as pro-objects in $\Fin\Set$ and as Stone spaces. The interplay between these two viewpoints is very useful for understanding monomorphisms and finite group actions, as the reader will see below. We denote the functor $\wh\Set \to \Set$ that sends a profinite set $\{X_i\}$ to its limit $\lim_i X_i$ by $U$; note that under the equivalence $\wh\Set \simeq \Stone$, this is the functor that sends a Stone space to its underlying set. Recall from \autoref{lemma:ReindexPosetSGA} that any pro-object is isomorphic to one indexed by a codirected poset.

\begin{lemma}\label{lem:EquivalentCharacterizationsMono}
Let $f \colon X \to Y$ be a map of profinite sets and assume that it is represented by a strict map $\{f_i \colon X_i \to Y_i\}_{i \in I}$ indexed by a codirected poset $I$. Then the following are equivalent:
\begin{enumerate}[(i)]
    \item \label{item1:EquivalentCharacterizationsMono} $f$ is a monomorphism in $\wh\Set$.
    \item \label{item2:EquivalentCharacterizationsMono}$\Delta \colon X \to X \times_Y X$ is an isomorphism.
    \item \label{item3:EquivalentCharacterizationsMono}For any $i \in I$, there exists a $j \leq i$ such that for any $x,x' \in X_j$, if $f_j(x) = f_j(x') \in Y_j$, then $\pi_{ij}(x) = \pi_{ij}(x')$ (and hence the same will hold for any $k \leq j$).
    \item \label{item4:EquivalentCharacterizationsMono}There exists a factorization $X \xrightarrow{\phi} X' \xrightarrow{g} Y$ of $f$, where $\phi$ and $g$ are again (represented by) strict maps, $\phi$ is an isomorphism, and $g$ is levelwise injective.
    \item \label{item5:EquivalentCharacterizationsMono}  $Uf \colon UX \to UY$ is an injective map of sets.
\end{enumerate}
\end{lemma}

\begin{proof}
The equivalence of \ref{item1:EquivalentCharacterizationsMono} and \ref{item2:EquivalentCharacterizationsMono} follows from general category theory. By writing out what it means for two maps in a pro-category to be equal, it follows that item \ref{item3:EquivalentCharacterizationsMono} is equivalent to the statement that $\Delta \circ \pi_1 \colon X \times_Y X \to X \times_Y X$ is the identity, which is clearly equivalent to \ref{item2:EquivalentCharacterizationsMono}. To see that \ref{item3:EquivalentCharacterizationsMono} implies \ref{item4:EquivalentCharacterizationsMono}, define $X'_i = \im(f_i) \subset Y_i$ and $X' = \{X'_i\}_{i \in I}$. Then $X \to X'$ is an isomorphism because of \ref{item3:EquivalentCharacterizationsMono}, and $X' \to Y$ is levelwise injective by construction. To see that \ref{item4:EquivalentCharacterizationsMono} implies \ref{item5:EquivalentCharacterizationsMono}, note that a cofiltered limit of monomorphisms is again a monomorphism in $\Set$. Finally, item \ref{item1:EquivalentCharacterizationsMono} follows from \ref{item5:EquivalentCharacterizationsMono} since a map in $\Stone$ is a monomorphism if and only if it is injective on the underlying sets.
\end{proof}

Throughout the rest of this section, let $G$ be a finite group. The equivalence $\wh\Set \simeq \Stone$ can easily be generalized to an equivalence between $\Pro(\Ac{G}{\Fin\Set})$ and $\GwhSet \simeq \Ac{G}{\Stone}$, where $\GwhSet$ and $\AcG{\Stone}$ are the categories of objects in $\wh\Set$ and $\Stone$, respectively, together with a (right) $G$-action. The equivalence $\Pro(\AcG{\Fin\Set}) \to \AcG{\Stone}$ is again defined by sending $\{X_i\}_{i \in I}$ to its limit $\lim_i X_i$ in the category of $G$-spaces, where the finite $G$-sets $X_i$ are endowed with the discrete topology. This can be seen as a consequence of \autoref{thm:MeyersResult}. In particular, if $G$ acts on a profinite set $X = \{X_i\}$, we may assume without loss of generality that this action is given by a levelwise action $G \times X_i \to X_i$. The objects of $\GwhSet$ will be called \emph{profinite $G$-sets}. We will say that $G$ acts \emph{freely} on a profinite set or Stone space $X$ if the map $X \times G \to X \times X; (x,g) \mapsto (x,xg)$ is a monomorphism.

\begin{lemma}\label{lem:EquivalentCharacterizationsFreeAction}
For a finite group $G$ acting levelwise on a profinite set $X = \{X_i\}_{i \in I}$ indexed by a codirected poset $I$, the following are equivalent:
\begin{enumerate}[(i)]
    \item \label{item1:EquivalentCharacterizationsFreeAction} The action is free.
    \item \label{item2:EquivalentCharacterizationsFreeAction} $G$ acts freely on the underlying set $UX = \lim_i X_i$.
    \item \label{item3:EquivalentCharacterizationsFreeAction} For each $i$, there exists a $j \leq i$ such that $G$ acts freely on $X_j$ (and hence also on each $X_k$ for $k \leq j$).
    \item \label{item4:EquivalentCharacterizationsFreeAction} $X$ is isomorphic to an object $X'$ with a levelwise free $G$-action.
\end{enumerate}
\end{lemma}

We leave the proof to the reader, with two remarks. First, the equivalence between \ref{item1:EquivalentCharacterizationsFreeAction} and \ref{item3:EquivalentCharacterizationsFreeAction} is special case of the equivalence between \ref{item1:EquivalentCharacterizationsMono} and \ref{item3:EquivalentCharacterizationsMono} of \autoref{lem:EquivalentCharacterizationsMono}. Secondly, item \ref{item4:EquivalentCharacterizationsFreeAction} follows from \ref{item3:EquivalentCharacterizationsFreeAction} by choosing any $j$ such that $G$ acts freely on $X_j$, and then defining $X' = \{X_k\}_{k \in I_{\leq j}}$, since the inclusion $I_{\leq j} \hookrightarrow I$ is a cofinal functor.

To study normal monomorphisms of dendroidal (profinite) sets, we will need a relative version of the above. In particular, for a map $X \to Y$ of profinite $G$-sets, we need to define what it means for $G$ to ``act freely on the complement''.

\begin{definition}\label{def:ActingFreelyOnComplement}
Let $f \colon X \to Y$ be a map between profinite $G$-sets. We say that $G$ \emph{acts freely on the complement} (of the image of $X$ in $Y$) if $G$ acts freely on $UY \setminus \im(Uf)$, where $Uf \colon UX \to UY$ is the map of underlying $G$-sets.
\end{definition}

This definition is the one that is the easiest to state, but in practice we will work with the equivalent characterizations given in below, partly because the complement of the image of $f$ occurring in the definition is not a Stone space. Since we will only consider the above definition for monomorphisms, we have added this extra assumption on $X \to Y$ in the following lemma.

\begin{lemma}\label{lem:EquivalentCharacterizationsFreeActionOnComplement}
Let $X \rightarrowtail Y$ be a monomorphism of profinite $G$-sets, and assume that it is represented by a strict map $\{f_i \colon X_i \to Y_i\}_{i \in I}$ between diagrams of finite $G$-sets indexed by a codirected poset $I$. Then the following are equivalent:
\begin{enumerate}[(i)]
    \item \label{item1:EquivalentCharacterizationsFreeActionOnComplement} $G$ acts freely on the complement.
    \item \label{item2:EquivalentCharacterizationsFreeActionOnComplement} For any non-unit $g \in G$, the equalizer $Y^g \rightarrowtail Y$ of the action by $g$ and the identity factors through $X$:
    \[\begin{tikzcd}
    Y^g \ar[r,tail] \ar[dr,dashed] & Y \ar[r, shift left = 0.4ex,"\id"] \ar[r,shift right = 0.4ex,"\cdot g"'] & Y \\
     & X \ar[u,tail] &
    \end{tikzcd}\]
    \item \label{item3:EquivalentCharacterizationsFreeActionOnComplement} For any $i \in I$, there exists a $j \leq i$ such that for any $g \in G$ and $y \in Y_j$, if $yg = y$, then either $g = e$ or $\pi_{ij}(y) \in f(X_i) \subset Y_i$ (and hence the same will hold for any $k \leq j$).
    \item \label{item4:EquivalentCharacterizationsFreeActionOnComplement} Up to isomorphism, $X \to Y$ admits a strict representation by injective maps of finite $G$-sets for which $G$ acts freely on the complement.
\end{enumerate}
\end{lemma}

\begin{proof}
The equivalence of \ref{item1:EquivalentCharacterizationsFreeActionOnComplement} and \ref{item2:EquivalentCharacterizationsFreeActionOnComplement} is a consequence of the fact that $U \colon \Stone \to \Set$ preserves equalizers. To see that \ref{item1:EquivalentCharacterizationsFreeActionOnComplement} implies \ref{item3:EquivalentCharacterizationsFreeActionOnComplement}, let $i \in I$ be given and let $y \in Y_i \setminus X_i$. Then $G$ acts freely on $\pi_i^{-1}(y) \subset UY$. Since $\pi_i^{-1}(y) = \lim_{j \leq i} \pi_{ij}^{-1}(y)$, we conclude by item \ref{item3:EquivalentCharacterizationsFreeAction} of \autoref{lem:EquivalentCharacterizationsFreeAction} that $G$ acts freely on $\pi_{ij}^{-1}(y)$ for some $j \leq i$ (and hence for any $k \leq j$). Since $Y_i$ is finite and $I$ is codirected, we can choose $j$ such that $G$ acts freely on $\pi_{ij}^{-1}(y)$ for any $y \in Y_i \setminus X_i$. It follows that, for any $g \in G$ and $y \in Y_j$ such that $yg = y$, either $g = e$ or $\pi_{ij}(y) \in f(X_i) \subset Y_i$.

We will now deduce \ref{item4:EquivalentCharacterizationsFreeActionOnComplement} from \ref{item3:EquivalentCharacterizationsFreeActionOnComplement}. By \autoref{lemma:ReindexPosetSGA}, we may assume that $\{f_i \colon X_i \to Y_i\}_i$ is indexed by a codirected poset $I$ such that $I_{\geq i}$ is finite for any $i \in I$. We may also assume without loss of generality that $f_i \colon X_i \to Y_i$ is injective for any $i \in I$ by replacing $X_i$ with its image in $Y_i$ under $f_i$ for any $i$, as in the proof of \ref{item4:EquivalentCharacterizationsMono} of \autoref{lem:EquivalentCharacterizationsMono}. We recursively define a monotone function $\theta \colon I \to I$ as follows. Let $i \in I$ be given and suppose that $\theta$ has been defined on the finite set $I_{\geq i}$. Choose any $j \leq i$ such that item \ref{item3:EquivalentCharacterizationsFreeActionOnComplement} of this lemma holds; note that since $\theta(I_{\geq i}) \subset I$ is finite, we can choose $j$ to be smaller than $\theta(k)$ for any $k \geq i$. Setting $\theta(i) = j$ defines a monotone and decreasing (hence cofinal) function $\theta \colon I \to I$, so we obtain canonical isomorphisms $\{X_{\theta(i)}\}_{i \in I} \cong X$ and $\{Y_{\theta(i)}\}_{i \in I} \cong Y$. For each $i \in I$, define $X'_i$ as the pullback
\[\begin{tikzcd}
X_{\theta(i)} \ar[drr,"f_{\theta(i)}", bend left = 13] \ar[ddr,"\pi_{i \theta(i)}"', bend right] \ar[dr, dashed,"\xi_i"] & & \\
& X'_i \arrow[dr,phantom, very near start, "\lrcorner"] \arrow[r,"\wt f_i", tail] \arrow[d] & Y_{\theta(i)} \arrow[d,"\pi_{i \theta(i)}"] \\
& X_i \arrow[r,"f_i", tail] & Y_i,
\end{tikzcd}\]
where the maps $\xi_i$ come from the universal property of the pullback. Then the maps $\xi_i \colon X_{\theta(i)} \to X'_i$ are natural in $i$, hence they define a levelwise map of profinite $G$-sets $\{X_{\theta(i)}\}_{i \in I} \to \{X'_i\}_{i \in I}$. The cofinality of $\theta \colon I \to I$ implies that the strict maps $\{X_{\theta(i)}\}_{i \in I} \to  \{X_i\}_{i \in I}$ and $\{Y_{\theta(i)}\}_{i \in I} \to \{Y_i\}_{i \in I}$ are isomorphisms of profinite $G$-sets. Since pullbacks in a pro-category can be computed ``levelwise'', the map $\{X'_i\}_{i \in I} \to \{X_i\}_{i \in I}$ is the pullback of an isomorphism of profinite $G$-sets, hence itself an isomorphism. By the 2 out of 3 property, we see that the strict map $\{X_{\theta(i)}\}_{i \in I} \to \{X'_i\}_{i \in I}$ defined by the maps $\xi_i$ is an isomorphism.
The map $\wt f_i$ is injective since it is a pullback of an injection. Furthermore, by construction of $\theta(i)$, it follows that $G$ acts freely on the complement of $\wt f_i \colon X'_i \to Y_{\theta(i)}$. In particular, $\{\wt f_i \colon X_i' \to Y_{\theta(i)}\}_{i \in I}$ is a strict map that is levelwise a monomorphism for which $G$ acts freely on the complement.

We leave it to the reader to show that \ref{item4:EquivalentCharacterizationsFreeActionOnComplement} implies \ref{item1:EquivalentCharacterizationsFreeActionOnComplement}.
\end{proof}

For the study of normal monomorphisms between dendroidal profinite sets in \autoref{sec:NormalMonos}, it will be useful to characterize those monomorphisms $f \colon X \to Y$ in $\GwhSet$ for which $G$ acts freely on the complement as those having the \llp with respect to a certain map. To this end, denote the two-point set by $\bftwo$ and denote the $G$-set obtained from the right action of $G$ on itself by $G$. Since $G$ is finite, we can view this both as an object of $\GSet$ and of $\GwhSet$. Recall that the forgetful functor $\AcG{\Fin\Set} \to \Fin\Set$ has a right adjoint $X \mapsto X^G$. (Note that this notation is not related to the $G$-fixed points, but that we use it as a shorthand for $\prod_{g \in G} X$.)

\begin{lemma}\label{lem:FreelyComplementAsLLP}
Let $f \colon X \to Y$ be a map of $G$-sets (resp. profinite $G$-sets). Then $f$ is a monomorphism for which $G$ acts freely on the complement if and only if $f$ has the \llp with respect to the coproduct of the two constant maps
\[\bftwo^G \sqcup G \to * \sqcup *\]
in $\GSet$ (resp. $\GwhSet$).
\end{lemma}

\begin{proof}
We leave the case where $f \colon X \to Y$ is a map of $G$-sets to the reader. Suppose now that $f \colon X \to Y$ is a map of profinite $G$-sets, and suppose that it has the \llp with respect to the map given in the statement of the lemma. Then $f \colon X \to Y$ in particular has the \llp with respect to $\bftwo^G \to *$, so the underlying map of profinite sets has the \llp with respect to $\bftwo \to *$. View $X$ and $Y$ as Stone spaces, and suppose that $f$ is not a monomorphism. Then there are $x,x' \in X$ such that $f(x) = f(x')$, yet $x \neq x'$. Choose some clopen $S$ around $x$ that does not contain $x'$ and let $\mathbb 1_S \colon X \to \bftwo$ denote the indicator function. This indicator function is continous, but a lift $h \colon Y \to \bftwo$ such that $\mathbb 1_S = hf$ clearly cannot exist, which is a contradiction. We conclude that $f \colon X \to Y$ is a monomorphism in $\Stone$, hence also a monomorphism in $\GwhSet$. To see that $G$ acts freely on the complement, again view $X$ and $Y$ as Stone spaces and let $y \in Y \setminus f(X)$ be given. Since $f(X)$ is closed in $Y$, there exists a clopen $S$ containing $y$ that is disjoint from $f(X)$. We then have a commutative square
\[\begin{tikzcd}
X \ar[r,"\alpha"] \ar[d,"f"] & \bftwo^G \sqcup G \ar[d] \\
Y \ar[r,"\mathbb 1_S"'] \ar[ur,dashed] & * \sqcup *,
\end{tikzcd}\]
where $\alpha$ sends $X$ to any point in the first summand of $\bftwo^G \sqcup G$. By assumption, a lift $h$ exists for this diagram. Since $h$ maps $y$ to a point in $G$, we see that $yg = y$ if and only if $g = e$.

For the converse, suppose that $f \colon X \to Y$ is a monomorphism for which $G$ acts freely on the complement. By item \ref{item4:EquivalentCharacterizationsFreeActionOnComplement} of \autoref{lem:EquivalentCharacterizationsFreeActionOnComplement}, we may assume that $f$ is represented by a strict map $\{f_i \colon X_i \cofarrow Y_i\}_{i \in I}$ such that, for every $i$, the map $f_i$ is an injection of finite $G$-sets for which $G$ acts freely on the complement. It follows from the case of (ordinary) $G$-sets that for any $i \in I$, the map $f_i$ has the \llp with respect to $\bftwo^G \sqcup G \to * \sqcup *$. Since finite $G$-sets are cocompact in $\GwhSet$, we conclude that $X \to Y$ must have the \llp with respect to this map as well.
\end{proof}

\section{Normal monomorphisms}\label{sec:NormalMonos}

The goal of this section is to define and study normal monomorphisms of dendroidal profinite sets. In particular, we will show that they form part of a (fibrantly generated) weak factorization system on $\dd\wh\Set$. They will play the role of cofibrations in the model structures constructed in \autoref{sec:TheModelStructure}. Unless stated otherwise, all results hold for open, closed and general dendroidal sets.

We will make use of the following fact regarding monomorphisms between dendroidal profinite sets.

\begin{lemma}\label{lem:MonoIsLevelwiseMonodSet}
	Let $f \colon X \to Y$ be a map of dendroidal profinite sets. Then $f$ is a monomorphism if and only if, up to isomorphism, it admits a level representation $\{f_i \colon X_i \to Y_i\}_{i \in I}$ such that for every $i \in I$, the map $f_i$ is a monomorphism between lean dendroidal sets.
\end{lemma}

\begin{proof}
	The ``if'' direction follows since cofiltered limits of monomorphisms are monomorphisms in $\wh\Set$. For the converse, let $f$ be a monomorphism and suppose without loss of generality that it has a level representation $\{f_i \colon X_i \to Y_i\}_{i \in I}$. As in the proof of \autoref{lem:EquivalentCharacterizationsMono}, if we replace $X_i$ by $X'_i = \im(f_i)$ for every $i \in I$, then $X$ is isomorphic to $\lim_i X_i'$ and $\{X_i' \hookrightarrow Y_i\}_{i \in I}$ is levelwise a monomorphism of degreewise finite dendroidal sets. Since $X' \cong \lim_n \cosk_n(X')$ and $Y = \lim_n \cosk_n(Y)$, we see that the diagram $\{\cosk_n(X_i') \hookrightarrow \cosk_n(Y_i)\}_{(i,n) \in I \times \bbN^{op}}$ is, up to isomorphism, a level representation of $f \colon X \to Y$ that is levelwise a monomorphism between lean dendroidal sets.
\end{proof}

A map $X \to Y$ of dendroidal sets is called a \emph{normal monomorphism} if it is a monomorphism for which $\Aut(T)$ acts freely on the complement of $X_T \hookrightarrow Y_T$, for each tree $T$. To generalize this to the profinite setting, recall from \autoref{def:ActingFreelyOnComplement} that if $f \colon X \to Y$ is a map of profinite $G$-sets, then $G$ is said to act freely on the complement if $G$ acts freely on the set $UY \setminus \im(Uf)$.

\begin{definition}
	A morphism $X \to Y$ of dendroidal profinite sets is called a \emph{normal monomorphism} if it is a monomorphism with the property that, for each tree $T$, the group $\Aut(T)$ acts freely on the complement of the image of $X_T \rightarrowtail Y_T$. A dendroidal profinite set is called \emph{normal} if the monomorphism $\varnothing \to Y$ is normal; i.e., if $\Aut(T)$ acts freely on $Y_T$ for every tree $T$.
\end{definition}

In order to prove that the normal monomorphisms of dendroidal profinite sets form part of a weak factorization system on $\dd\wh\Set$, it would be useful to have a characterization of them similar to \autoref{lem:MonoIsLevelwiseMonodSet}; that is, we would like that $f \colon X \to Y$ is a normal monomorphism if and only if up to isomorphism, it admits a strict representation $\{f_i \colon X_i \to Y_i\}$ by normal monomorphisms of lean dendroidal sets. However, the following example shows that this notion is too restrictive when working with lean dendroidal sets.

\begin{example}\label{exx:LeanDendroidalSetIsNotNormal}
We claim that a monomorphism between lean dendroidal sets is normal if and only if it is an isomorphism. In particular, by considering the map $\varnothing \to X$, we see that a lean dendroidal set is normal if and only if it is the initial object. To see that this holds, let $f \colon X \to Y$ be a normal monomorphism between lean dendroidal sets and suppose that for some tree $T$, there exists an element $y \in Y_T$ that is not in the image of $f$. Let $n \geq 2$ be such that $|T| \leq n$ and such that $Y$ is $n$-coskeletal. Construct the tree $C_{n,T}$ by grafting $T$ onto the $(n+1)$-st leaf of the corolla $C_{n+1}$; this tree is pictured below.
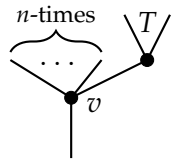
\[\begin{tikzpicture}
    \draw[thick] (0,0) -- (0,0.8);
    \draw[thick] (0,0.8) -- (-0.8,1.3);
    \draw[thick] (0,0.8) -- (0.4,1.3);
    \draw[thick] (0,0.8) -- (1,1.3);
    \draw[thick] (1,1.3) -- (0.7,1.9);
    \draw[thick] (1,1.3) -- (1.3,1.9);
    
    \draw[fill] (0,0.8) circle [radius=0.085];
    \draw[fill] (1,1.3) circle [radius=0.085];
    
    \node at (-0.15,1.25) {$\cdots$};
    \node at (1,1.8) {$T$};
    \node at (0.3,0.7) {$v$};
    
    \draw [thick,decorate,decoration={brace,amplitude=8pt},xshift=0pt,yshift=3pt](-0.8,1.3) -- (0.4,1.3) node[black,midway,yshift=0.5cm] {\footnotesize $n$-times};
\end{tikzpicture}\]
Denote the root vertex of this tree by $v$. The root face of this tree is the inclusion $\delta_v \colon T \hookrightarrow C_{n,T}$. Note that $\Sigma_n \hookrightarrow \Aut(C_{n,T})$, where $\Sigma_n$ acts by permuting the $n$ leaves attached to the root vertex. Since the root vertex of $C_{n,T}$ has $n+1$ incoming edges, we see for any tree $S$ that if there exists a map $S \to C_{n,T}$ that sends the root edge of $S$ to the root edge of $C_{n,T}$, then either $S$ factors through the root inclusion $\eta \hookrightarrow C_{n,T}$ or $S$ has at least $n$ leaves. This implies that any map $S \to C_{n,T}$ with $|S| \leq n$ either factors through $\delta_v \colon T \hookrightarrow C_{n,T}$ or through the inclusion $\eta \hookrightarrow C_{n,T}$ of the root edge or one of the leaves attached to the root vertex. In particular, these maps form a (discrete) final subcategory of $\Omega_{(n)}/C_{n,T}$. It follows that $Y_{C_{n,T}} \cong Y_\eta \times (Y_\eta)^n \times Y_T$ and that $\Sigma_n \subset \Aut(C_{n,T})$ acts by permuting the $n$ copies of $Y_\eta$ corresponding to the leaves. Now let $c \in Y_\eta$ be any colour, e.g. the colour corresponding to the root of $y \in Y_T$. Then the tuple $(c,\ldots,c,y)$ defines an element of $Y_{C_{n,T}} \cong Y_\eta \times (Y_\eta)^n \times Y_T$ that is not in the image of $f \colon X_{C_{n,T}} \to Y_{C_{n,T}}$. However, the subgroup $\Sigma_n \subset \Aut(C_{n,T})$ fixes this element, contradicting the fact that $f$ is normal. We conclude that $X_T \to Y_T$ must be surjective for any tree $T$, and hence that $f$ is an isomorphism.
\end{example}

\begin{remark}
	Note that the preceding argument works for open and general dendroidal sets, but not for closed dendroidal sets. In fact, there do exist normal monomorphisms between lean closed dendroidal sets that are not isomorphisms, as illustrated by the fact that the closed nerve of the associative operad is a normal lean closed dendroidal set (cf. \autoref{exx:AssociativeOperadCoskeletal}).
\end{remark}

This motivates the following definition.

\begin{definition}
\begin{enumerate}[(a)]
    \item A dendroidal set $X$ is called \emph{$n$-partially normal}, or \emph{$n$-normal} for short, if $\Aut(T)$ acts freely on $X_T$ for each tree $T$ of size at most $n$.
    \item Similarly, a morphism $f \colon X \to Y$ of dendroidal sets is called \emph{$n$-partially normal} or \emph{$n$-normal} if, for each tree $T$ of size at most $n$, the map $X_T \to Y_T$ is a monomorphism and $\Aut(T)$ acts freely on $Y_T \setminus f(X_T)$.
\end{enumerate}
\end{definition}

Partially normal morphisms can equivalently be defined in terms of skeleta.

\begin{lemma}\label{lem:nSkeletonOfNNormalMapIsNormal}
A morphism $X \to Y$ is $n$-normal if and only if $\sk_n(X) \to \sk_n(Y)$ is a normal monomorphism. In particular, a dendroidal set $X$ is $n$-normal if and only if $\sk_n(X)$ is normal.
\end{lemma}

\begin{proof}
Recall from \autoref{ssec:LeanDendroidalSets} that the inclusion $\sk_n(X) \hookrightarrow X$ induces an isomorphism $\sk_n (X)_T \to X_T$ for any tree $T$ with $|T| \leq n$. In particular, if $\sk_n(X) \to \sk_n(Y)$ is a normal monomorphism, then $X \to Y$ is $n$-normal. Furthermore, for the converse it suffices to prove that if $X \to Y$ is $n$-normal, then for any tree $T$ with $|T| > n$, the map $\sk_n(X)_T \to \sk_n(Y)_T$ is a monomorphisms for which $\Aut(T)$ acts freely on the complement.

To see that this map is a monomorphism, note that as mentioned in \autoref{ssec:LeanDendroidalSets}, $\sk_n X$ is a subobject of $X$.

Now note that for any tree $T \in \Omega$, the group $\Aut(T)$ acts by precomposition on the set of degeneracies $T \twoheadrightarrow S$. It follows from Proposition 3.13 of \cite{HeutsMoerdijk2020Trees} that this action is free. Combining this with the fact that up to isomorphism, there is an essentially unique way to write an element of $\sk_n(Y)_T$ as $\sigma^*(y)$ with $\sigma \colon T \twoheadrightarrow S$ a degeneracy and $y \in \sk_n(Y)_S = Y_S$ non-degenerate, we see that $\Aut(T)$ must act freely on the complement of $\sk_n(X)_T \to \sk_n(X)_T$ if $X \to Y$ is $n$-normal and $|T| > n$.
\end{proof}

While we cannot represent normal monomorphisms by levelwise normal monomorphisms of lean dendroidal sets, the next best thing is possible: any normal monomorphism admits, up to isomorphism, an \emph{increasingly normal representation}.

\begin{definition}
	Let $X \to Y$ be a map of dendroidal profinite sets and assume that it is represented by a strict map $\{f_i \colon X_i \to Y_i\}_{i \in I}$ between diagrams of lean dendroidal sets indexed by a codirected poset $I$. We say that $\{f_i \colon X_i \to Y_i\}$ is an \emph{increasingly normal representation} if for any $n \in \bbN$, there exists an $i \in I$ such that for any $j \leq i$, the map $f_j \colon X_j \to Y_j$ is $n$-partially normal.
\end{definition}

\begin{lemma}\label{lem:NormalMonoHasIncreasinglyNormalRepresentation}
	Let $f \colon X \to Y$ be a map of dendroidal profinite sets. Then $f$ is a normal monomorphism if and only if, up to isomorphism, it admits an increasingly normal representation.
\end{lemma}

\begin{proof}
	For the ``if'' direction, suppose that $\{f_i \colon X_i \to Y_i\}$ is an increasingly normal representation of a map $f \colon X \to Y$ of dendroidal sets. Then for any tree $T$, there exists an $i$ such that for any $j \leq i$, the map $(X_j)_T \to (Y_j)_T$ is a monomorphism for which $\Aut(T)$ acts freely on the complement. Since $I_{\leq i} \subset I$ is cofinal, we see by \autoref{lem:EquivalentCharacterizationsFreeActionOnComplement} that $X_T \to Y_T$ is a monomorphism for which $\Aut(T)$ acts freely on the complement. We conclude that $X \to Y$ is normal.
	
	The proof of the converse direction is similar to that of item \ref{item4:EquivalentCharacterizationsFreeActionOnComplement} of \autoref{lem:EquivalentCharacterizationsFreeActionOnComplement}, with a few subtle differences. Let $f \colon X \to Y$ be a normal monomorphism of dendroidal profinite sets. By \autoref{lem:MonoIsLevelwiseMonodSet} and \autoref{lemma:ReindexPosetSGA}, we may assume that it admits a level representation $\{f_i \colon X_i \hookrightarrow Y_i\}_{i \in I}$ by injective maps indexed by a codirected poset $I$ such that for every $i \in I$, the set $I_{> i}$ is finite. In particular, we can define a strictly order-reversing map $\phi \colon I \to \bbN$ by setting $\phi(i) = \# I_{>i}$. We recursively define a decreasing monotone function $\theta \colon I \to I$. Let $i \in I$ be given and suppose that $\theta$ has been defined on $I_{> i}$. Define $\theta(i)$ to be any element in $I$ such that $\theta(i) \leq i$, such that $\theta(i) < \theta(k)$ for all $k \in I_{>i}$, and such that for any tree $T$ of size $\leq \phi(i)$, item \ref{item3:EquivalentCharacterizationsFreeActionOnComplement} of \autoref{lem:EquivalentCharacterizationsFreeActionOnComplement} holds for the map of $\Aut(T)$-sets $(X_i)_T \hookrightarrow (Y_i)_T$ with $j = \theta(i)$. Such an element $\theta(i)$ exists by codirectedness of $I$. For each $i \in I$, define $X'_i$ as the pullback
	\[\begin{tikzcd}
		X'_i \arrow[dr,phantom, very near start, "\lrcorner"] \arrow[r, hook] \arrow[d] & Y_{\theta(i)} \arrow[d,"\pi_{i \theta(i)}"] \\
		X_i \arrow[r,"f_i", hook] & Y_i.
	\end{tikzcd}\]
	Exactly as in the proof item \ref{item4:EquivalentCharacterizationsFreeActionOnComplement} of \autoref{lem:EquivalentCharacterizationsFreeActionOnComplement}, it follows that up to isomorphism, $\{X'_i \hookrightarrow Y_{\theta(i)}\}_{i \in I}$ is a level representation of $f \colon X \to Y$. Furthermore, by construction the action of $\Aut(T)$ on the complement of $(X_i')_T \hookrightarrow (Y_{\theta(i)})_T$ is free whenever $|T| \leq \phi(i)$, so $X'_i \to Y_{\theta(i)}$ is $n$-normal for $n = \phi(i)$. We conclude that $\{X_i' \hookrightarrow Y_{\theta(i)}\}_{i \in I}$ is the desired increasingly normal representation.
\end{proof}

This lemma can be used to characterize the normal monomorphisms in $\dd\wh\Set$ as those maps having the \llp \wrt certain maps. Recall from \autoref{ssec:dSetsAndSpaces} that the fibrant objects in the operadic model structure on $\dd\Set$ are called $\infty$-operads.

\begin{proposition}\label{prop:NormalMonoIffLLPwrtLeanTrivFibs}
	A morphism $X \to Y$ of dendroidal profinite sets is a normal monomorphism if and only if it has the \llp \wrt all trivial fibrations between lean $\infty$-operads.
\end{proposition}

\begin{proof}
	For a (profinite) right $\Aut(T)$-set $Z$, viewed as a functor out of $\Aut(T)^{op}$, denote the right Kan extension along the inclusion $\Aut(T)^{op} \hookrightarrow \Omega^{op}$ by $R_T Z$. Recall the map $\bftwo^G \sqcup G \to * \sqcup *$ from \autoref{lem:FreelyComplementAsLLP}. We denote the map $R_T(\bftwo^{\Aut(T)} \sqcup \Aut(T)) \to R_T(* \sqcup *)$ by $\Psi_T$. It follows from \autoref{lem:FreelyComplementAsLLP} that a map of dendroidal (profinite) sets is a normal monomorphism if and only if it has the \llp \wrt all maps in the set $\{\Psi_T\}_{T \in \Omega}$. In particular, the maps $\Psi_T$ are trivial fibrations. It is easy to show that the domain of $R_T(* \sqcup *)$ is a (weakly contractible) lean $\infty$-operad, so we conclude that if a map of dendroidal profinite sets $X \to Y$ has the \llp \wrt trivial fibrations between lean $\infty$-operads, then it is a normal monomorphism.
	
	For the converse, suppose that $X \cofarrow Y$ is normal and that $B \trivfibarrow A$ is a trivial fibration between lean $\infty$-operads. By \autoref{lem:NormalMonoHasIncreasinglyNormalRepresentation}, we may assume that $X \cofarrow Y$ admits an increasingly normal representation $\{X_i \to Y_i\}$. Since $B$ and $A$ are cocompact, it follows that $X \cofarrow Y$ has the \llp \wrt $B \trivfibarrow A$ if $X_i \to Y_i$ does for small enough $i$. Since there exists an $n$ such that $A$ and $B$ are $n$-coskeletal, it suffices by adjunction to show that we can construct a lift in
	\[\begin{tikzcd}
		\sk_n(X_i) \ar[r] \ar[d] & B \ar[d,two heads,"\sim" rot90] \\
		\sk_n(Y_i) \ar[r] & A.
	\end{tikzcd}\]
	for small enough $i$. Since $\{X_i \to Y_i\}$ is increasingly normal, it follows from \autoref{lem:nSkeletonOfNNormalMapIsNormal} that the left-hand vertical map is a normal monomorphism for small enough $i$, hence the operadic model structure on $\dd\Set$ provides the desired lift.
\end{proof}

\begin{remark}\label{remark:NormalMonoIffLLPOtherSets}
	The proof of \autoref{prop:NormalMonoIffLLPwrtLeanTrivFibs} actually shows something stronger than the statement of that proposition: a map of dendroidal profinite sets is already a normal monomorphism if it has the \llp \wrt trivial fibrations between weakly contractible lean $\infty$-operads, and a normal monomorphism of dendroidal profinite sets has the \llp \wrt every trivial fibration between lean dendroidal sets.
\end{remark}

\section{The model structures for profinite \texorpdfstring{$\infty$}{infinity}-operads} \label{sec:TheModelStructure}

We start this section by constructing a convenient normalization functor for dendroidal profinite sets. This functor is then used to define (operadic) weak equivalences of dendroidal profinite sets and prove some of their basic properties. We conclude this section by showing that these are the weak equivalences of a model structure on $\dd\wh\Set$ in which the cofibrations are the normal monomorphisms.

For simplicity of exposition, most results in this section are only stated for general dendroidal sets. However, unless stated otherwise, these results also hold for open and closed dendroidal sets.

\subsection{A convenient normalization functor} \label{ssec:Normalization}

Recall from \autoref{ssec:dSetsAndSpaces} that for dendroidal sets $X$ and $Y$, one defines the simplicial hom $\sHom(X,Y)$ by
\[\sHom(X,Y)_\bullet = \Hom(X \otimes \Delta[\bullet], Y) \cong \Hom(X,Y^{\Delta[\bullet]}), \]
where $\Delta[\bullet]$ is identified with $i_! \Delta[\bullet]$ and $\Hom$ is the set of morphisms in $\dd\Set$. A map $X \to Y$ of dendroidal sets is called an operadic weak equivalence if for any $\infty$-operad $A$, the map
\[\sHom(\wt Y, A) \to \sHom(\wt X, A) \]
is a weak equivalence in the Joyal model structure on $\s\Set$, where $\wt X$ and $\wt Y$ denote (functorial) normalizations of $X$ and $Y$. By a normalization of a dendroidal set $X$, we mean a normal dendroidal set $\wt X$ together with a trivial fibration $\wt X \to X$. In order to mimic this definition for dendroidal profinite sets, we need a functorial normalization in $\dd\wh\Set$. In light of \autoref{prop:NormalMonoIffLLPwrtLeanTrivFibs}, we could use the cosmall object argument to construct such a functorial normalization. However, since we have very little control over the normalizations obtained in this manner, we use a different construction based on the observation that the product of a dendroidal set with a normal dendroidal set is always normal.

\begin{lemma}\label{lemma:DegreewiseFiniteNormalizationPoint}
There exists a degreewise finite normal dendroidal set $E$ such that the map  $E \to *$ to the terminal object is a trivial fibration (i.e. there exists a weakly contractible degreewise finite normal $\infty$-operad).
\end{lemma}

\begin{remark}
	For open and general dendroidal sets, one can use the homotopy-coherent nerve of the (open or unital) Barratt-Eccles operad as a degreewise finite normalization of the point. However, since the unital Barratt-Eccles operad is not Reedy fibrant, its closed homotopy-coherent nerve is not a closed $\infty$-operad (i.e. a fibrant object in $\cd\Set$). For this reason we use the more abstract proof of \autoref{lemma:DegreewiseFiniteNormalizationPoint} given below.
\end{remark}

\begin{proof}[Proof of \autoref{lemma:DegreewiseFiniteNormalizationPoint}]
	We construct the desired dendroidal set $E$ by modifying Quillen's small object argument. Start by defining $E^{(-1)} = \varnothing$. Let $n \in \bbN$ and suppose that $E^{(n-1)}$ has been constructed. Define $E^{(n)}$ as the pushout
	\begin{equation}\label{diag:PushoutSmallObjects}
		\begin{tikzcd}
			\coprod \partial \Omega[T] \ar[r] \ar[d] & E^{(n-1)} \ar[d] \\
			\coprod \Omega[T] \ar[r] & E^{(n)}
		\end{tikzcd}
	\end{equation}
	where the coproduct ranges over all maps $\partial \Omega[T] \to E^{(n-1)}$ with the property that $|T| = n$. (Note that we do not take all maps $\partial \Omega[T] \to E^{(n-1)}$ for \emph{all} trees $T$, as in the usual small object argument.) Define $E = \colim_{n \in \bbN} E^{(n)}$. Then $E$ is normal since $\partial \Omega[T] \to \Omega[T]$ is a normal monomorphism for any tree $T$.
	
	To show that $E \to *$ is a trivial fibration, recall from \autoref{ssec:LeanDendroidalSets} that $\sk_{|T|-1} \Omega[T] = \partial \Omega[T]$. In particular, the inclusion $\partial \Omega[T] \to \Omega[T]$ induces an isomorphism $(\partial \Omega[T])_S \cong \Omega[T]_S$ for any tree $S$ with $|S| < |T|$. Since the pushout \eqref{diag:PushoutSmallObjects} is computed levelwise, this implies that for any tree $S$ of size $\leq n$, the map $E^{(n)}_S \to E^{(n+1)}_S$ is an isomorphism. In particular, $E^{(n)} \to E$ induces an isomorphism on $n$-skeleta. Now suppose that a map $\partial \Omega[T] \to E$ is given, and let $n = |T|$. Since $\partial \Omega[T] = \sk_{n-1}\Omega[T]$ is $(n-1)$-skeletal, we see that this map must factor through $E^{(n-1)}$. By construction of the pushout \eqref{diag:PushoutSmallObjects}, the map $\partial \Omega[T] \to E^{(n-1)}$ admits a filler in $E^{(n)}$, hence the original map $\partial \Omega[T] \to E$ admits a filler. We conclude that $E \to *$ is a trivial fibration.
	
	It follows by induction that the coproduct in the pushout \eqref{diag:PushoutSmallObjects} is always finite and hence that $E^{(n)}$ is degreewise finite and skeletal for every $n \in \bbN$. In particular, since $E^{(n)} \to E$ is an isomorphism on $n$-skeleta, it follows that $E$ is degreewise finite.
\end{proof}

Throughout the rest of this section, let $E$ denote the degreewise finite normalization of the terminal object $*$ constructed in \autoref{lemma:DegreewiseFiniteNormalizationPoint}. 

\begin{lemma}
For any dendroidal profinite set $X$, the product $X \times E$ is normal and the projection $X \times E \to X$ has the \rlp \wrt any normal monomorphism in $\dd\wh\Set$.
\end{lemma}

\begin{proof}
As $X \times E$ is clearly normal, it suffices to show that $E \to *$ has the \rlp \wrt any normal monomorphism in $\dd\wh\Set$. As described in \autoref{ssec:LeanDendroidalSets}, $E$ is the limit of its tower of coskeleta
\[\cosk_0 E \leftarrow \cosk_1 E \leftarrow \cosk_2 E \leftarrow \ldots \leftarrow E\]
so it suffices to show that $\cosk_n E \to \cosk_{n-1} E$ has the \rlp \wrt any normal monomorphism in $\dd\wh\Set$. Since this is a map between lean dendroidal sets, it suffices by \autoref{prop:NormalMonoIffLLPwrtLeanTrivFibs} to show that it has the \rlp \wrt $\partial \Omega[T] \to \Omega[T]$ for every tree $T$. By adjunction, this is equivalent to $E$ having the \rlp \wrt
\[\sk_n \partial \Omega[T] \cup_{\sk_{n-1} \partial \Omega[T]} \sk_{n-1} \Omega[T] \to \sk_n \Omega[T]\]
which is an isomorphism if $|T| \neq n$ and the inclusion $\partial \Omega[T] \hookrightarrow \Omega[T]$ if $n =|T|$.
\end{proof}

\subsection{Weak equivalences of dendroidal profinite sets} \label{ssec:WeakEquivalencesDSets}

The simplicial hom defined in \autoref{ssec:Normalization} is part of a two-variable adjunction, meaning that there exist tensor and cotensor functors
\[\otimes \colon \dd\Set \times \s\Set \to \dd\Set \quad \text{and} \quad (\mhyphen)^{(\mhyphen)} \colon \s\Set^{op} \times \dd\Set \to \dd\Set  \]
together with natural isomorphisms
\[\Hom(X \otimes M, Y) \cong \Hom(M,\sHom(X,Y)) \cong \Hom(X, Y^M).\]
The tensor and cotensor are (by slight abuse of notation) defined by $X \otimes M := X \otimes i_! M$ and $X^M := X^{i_! M}$, respectively, where $i_!$ denotes the inclusion $\s\Set \hookrightarrow \dd\Set$.

By \autoref{lemma:TrueExponentialLean}, the cotensor restricts to a functor $\FinsSet^{op} \times \LdSet \to \LdSet \hookrightarrow \dd\wh\Set$ that preserves finite limits in both variables. Here $\FinsSet$ denotes the category of finite simplicial sets; that is, simplicial sets that have finitely many non-degenerate simplices. Since any simplicial set is the union of its finite simplicial subsets and since finite simplicial sets are compact in $\s\Set$, we see that $\Pro(\FinsSet^{op}) \simeq \s\Set^{op}$. As explained in \autoref{ssec:ProCategories}, this implies that the cotensor extends to a functor $(\mhyphen)^{(\mhyphen)} \colon \s\Set^{op} \times \dd\wh\Set \to \dd\wh\Set$ that is part of a two-variable adjunction, meaning that there exist functors
\[\otimes \colon \dd\wh\Set \times \s\Set \to \dd\wh\Set \quad \text{and} \quad \sHom(-,-) \colon \dd\wh\Set^{op} \times \dd\wh\Set \to \s\Set  \]
together with natural isomorphisms
\begin{equation*}
\Hom(X \otimes M, Y) \cong \Hom(M,\sHom(X,Y)) \cong \Hom(X, Y^M).
\end{equation*}
Using these isomorphisms and the Yoneda lemma, one deduces the following description of the simplicial hom of $\dd\wh\Set$ for $X = \{X_i\}$ and $Y = \{Y_j\}$:
\[ \sHom(X,Y) \cong \lim_j \colim_i \sHom(X_i, Y_j). \]
Here $\sHom$ on the right-hand side denotes the usual simplicial hom of $\dd\Set$, restricted to the full subcategory $\LdSet$. In particular, by definition of the morphisms in a pro-category, we obtain isomorphisms $\sHom(X,Y)_0 \cong \Hom(X,Y)$ that are natural in $X$ and $Y$.

We will now prove a weak version of the pullback-power property for this simplicial hom. To clarify the following lemma, let $H$ denote the simplicial set obtained by gluing two $2$-simplices to each other along the edges opposite to the $0$th and $2$nd vertex, respectively, and then collapsing the edges opposite to the $1$st vertex to a point in both of these $2$-simplices. This means that $H$ looks as follows, where the dashed lines represent the collapsed edges:
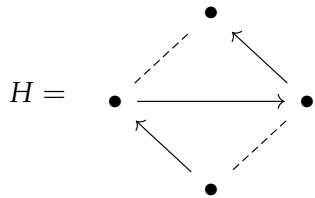
\[\begin{array}{rl}
   H =  &  \begin{tikzcd}[column sep = scriptsize]
 & \bullet & \\
\bullet \ar[ur, dashed, no head] \ar[rr] & & \bullet \ar[ul] \ar[dl, dashed, no head] \\
 & \bullet \ar[ul] &
\end{tikzcd}
\end{array} \]
By an \emph{$H$-fibration}, we mean an inner fibration of simplicial sets that has the \rlp \wrt $\{0\} \hookrightarrow H$. It follows from \cite[Proposition 2.4.6.5]{Lurie2009HTT} that for any $\infty$-category $C$, a map of simplicial sets $D \to C$ is an $H$-fibration if and only if it is a categorical fibration (cf. \cite[Lemma 2.1]{BlomMoerdijk2020SimplicialProV1}).

\begin{lemma}\label{lemma:PullbackPowerIsJFibration}
Let $X \to Y$ be a normal monomorphism in $\dd\wh\Set$ and let $B \fibarrow A$ be an operadic fibration between lean $\infty$-operads. Then the map
\begin{equation}\label{eq:PullbackPower0}\sHom(Y,B) \to \sHom(X,B) \times_{\sHom(X,A)} \sHom(Y,A) \end{equation}
is an $H$-fibration of simplicial sets, which is a trivial fibration if $B \fibarrow A$ is trivial or if $X \to Y$ has the \llp \wrt any operadic fibration between lean $\infty$-operads.
\end{lemma}

\begin{proof}
Since the inner horns $\Lambda^i[n] \to \Delta[n]$ and the endpoint inclusion $\{0\} \hookrightarrow H$ are trivial cofibrations in the Joyal model structure, we see that for any operadic fibration between lean $\infty$-operads $B \fibarrow A$, the pullback-power maps $B^{\Delta[n]} \fibarrow B^{\Lambda^i[n]} \times_{A^{\Lambda^i[n]}} A^{\Delta[n]}$ and $B^H \fibarrow B \times_A A^H$ are trivial fibrations between $\infty$-operads. Furthermore, the domains and codomains of these maps are lean by \autoref{lemma:TrueExponentialLean}. In particular, if $X \to Y$ is a normal monomorphism in $\dd\wh\Set$, then it follows by adjunction from \autoref{prop:NormalMonoIffLLPwrtLeanTrivFibs} that the map \eqref{eq:PullbackPower0} is an $H$-fibration.

Now suppose that $B \fibarrow A$ is a trivial fibration or that $X \to Y$ has the \llp \wrt any operadic fibration between lean $\infty$-operads. In this case an argument similar to the one above, but with the boundary inclusions $\partial \Delta[n] \to \Delta[n]$ in place of the inner horns and $\{0\} \hookrightarrow H$, shows that the map \eqref{eq:PullbackPower0} is a trivial fibration.
\end{proof}

\begin{corollary}\label{cor:PullbackPowerBetweenNormal}
Let $X \to Y$ be a normal monomorphism between normal dendroidal profinite sets and let $B \fibarrow A$ be an operadic fibration between lean $\infty$-operads. Then the map \eqref{eq:PullbackPower0} is a categorical fibration between $\infty$-categories.
\end{corollary}

\begin{proof}
As mentioned above, if $C$ is an $\infty$-category, then a map $D \to C$ from another simplicial set $D$ is an $H$-fibration if and only if it is a categorical fibration. Several applications of \autoref{lemma:PullbackPowerIsJFibration} show that the map \eqref{eq:PullbackPower0} is an $H$-fibration and that its codomain is an $\infty$-category, hence the result follows.
\end{proof}

\begin{remark}
	In \autoref{prop:PullbackPowerProperty} below, we show that \eqref{eq:PullbackPower0} is a categorical fibration under much more general hypotheses.
\end{remark}

Recall that $E$ denotes a degreewise finite normalization of the terminal object.

\begin{definition}
A map $X \to Y$ in $\dd\wh\Set$ is a called a \emph{weak equivalence} if for any lean $\infty$-operad $A$, the induced map
\[\sHom(Y \times E, A) \to \sHom(X \times E, A) \]
is a weak equivalence of $\infty$-categories.
\end{definition}

\begin{remark}
	It would also be reasonable to define a map $X\to Y$ to be a weak equivalence if and only if for any lean $\infty$-operad $A$, the induced map
	\[\Map(Y \times E, A) \to \Map(X \times E, A) \]
	is a weak equivalence of Kan complexes. Here $\Map(- \times E ,A)$ denotes the maximal Kan complex $k \sHom(- \times E, A)$ contained in the $\infty$-category $\sHom(- \times E, Z)$. It can be shown that this notion of weak equivalence is equivalent to the one defined above.
\end{remark}

We will see below that this definition does not depend on the choice of the degreewise finite normalization $E$ (cf. \autoref{prop:AlternativeNormalizations}). However, until then we let $E$ be fixed.

\begin{proposition}\label{prop:WeakEquivalencesStableUnderCofilLimits}
Weak equivalences in $\dd\wh\Set$ are stable under cofiltered limits.
\end{proposition}

\begin{proof}
Let $\{X_i \wearrow Y_i\}_{i \in I}$ be a natural weak equivalence between diagrams of dendroidal profinite sets indexed by a cofiltered category $I$, and write $X = \lim_i X_i$ and $Y = \lim_i Y_i$. Let $A$ be a lean $\infty$-operad. Since $A$ is cocompact in $\dd\wh\Set$, the horizontal maps in
\[\begin{tikzcd}
\colim\limits_{i \in I} \sHom(Y_i \times E, A) \ar[r,"\cong"] \ar[d] & \sHom(Y \times E, A) \ar[d] \\
\colim\limits_{i \in I} \sHom(X_i \times E, A) \ar[r, "\cong"] & \sHom(X \times E, A)
\end{tikzcd}\]
are isomorphisms. Since weak equivalences are stable under filtered colimits in the Joyal model structure on $\s\Set$, we see that the left-hand vertical map is a weak equivalence, hence $\sHom(Y \times E, A) \to \sHom(X \times E, A)$ is a weak equivalence. We conclude that $X \to Y$ is a weak equivalence of dendroidal profinite sets.
\end{proof}

\begin{remark}
It follows from the definition of the weak equivalences in $\dd\Set$ that the embedding $\LdSet \hookrightarrow \dd\wh\Set$ preserves weak equivalences. In particular, the previous proposition implies that if a morphism $X \to Y$ in $\dd\wh\Set$ has a strict representation $\{f_i \colon X_i \to Y_i\}$ for which each $X_i \to Y_i$ is an operadic weak equivalence of lean dendroidal sets, then $X \to Y$ is a weak equivalence in $\dd\wh\Set$.
\end{remark}

\begin{proposition}\label{prop:NormalMonoOperadicEquivalenceIfLLPwrtFibrations}
Let $X \to Y$ be a map in $\dd\wh\Set$. If $X \to Y$ has the \llp \wrt all fibrations between lean $\infty$-operads, then it is a normal monomorphism and a weak equivalence.
\end{proposition}

\begin{proof}
	It follows from \autoref{prop:NormalMonoIffLLPwrtLeanTrivFibs} that $X \to Y$ is a normal monomorphism. We claim that for any lean $\infty$-operad $A$, the map $\sHom(Y \times E,A) \to \sHom(X \times E, A)$ is a trivial fibration of simplicial sets, so in particular a weak equivalence. Note that this map has the \rlp \wrt $\partial \Delta[n] \to \Delta[n]$ if and only if $X \times E \to Y \times E$ has the \llp \wrt $A^{\Delta[n]} \to A^{\partial \Delta[n]}$. Since the latter map is an operadic fibration between lean $\infty$-operads, it suffices to show that if $X \to Y$ has the \llp \wrt all fibrations between lean $\infty$-operads, then so does $X \times E \to Y \times E$. By adjunction, this comes down to showing that if $B \fibarrow A$ is an operadic fibration between lean $\infty$-operads, then $\CExp(E,B) \to \CExp(E,A)$ is as well. It follows from \autoref{lemma:CartesianExponentialLean} that this is a map between lean dendroidal sets, so it suffices to show that $\CExp(E,-)$ preserves fibrations. This in turn follows from the fact that $- \times E$ preserves trivial cofibrations in the operadic model structure $\dd\Set$.
\end{proof}

\subsection{Construction of the model structures}\label{ssec:ModelStructureDSets}

We are now in a position to derive the existence of a Quillen model structure on the category $\dd\wh\Set$ of dendroidal profinite sets and its open and closed analogues, with the weak equivalences as defined in the previous section. We formulate this explicitly as follows.

\begin{theorem}\label{theorem:ProfiniteOperadicModelStructure}
There exist fibrantly generated model structures on $\dd\wh\Set$, $\od\wh\Set$ and $\cd\wh\Set$ in which the cofibrations are the normal monomorphisms and a map $X \to Y$ is a weak equivalence if and only if for any lean (open/closed) $\infty$-operad $A$, the map
\[\sHom(Y \times E, A) \to \sHom(X \times E, A) \]
is a weak equivalence of $\infty$-categories. These model structures are left proper.
\end{theorem}

The model structure on $\dd\wh\Set$ described in this theorem will be called the \emph{model structure for profinite $\infty$-operads}. Similarly, the model structures on $\od\wh\Set$ and $\cd\wh\Set$ will be called the \emph{model structures for open and closed profinite $\infty$-operads}, respectively. We will also call these model structures the \emph{operadic model structures} (on $\dd\wh\Set$, $\od\wh\Set$ and $\cd\wh\Set$).

This theorem is proved by choosing specific sets of generating (trivial) fibrations and checking the hypotheses of (the dual of) Kan's recognition theorem, spelled out as Theorem 11.3.1 in \cite{Hirschhorn2003Model}. For general dendroidal sets, the set of generating fibrations is given by
\[\mathcal{P} = \{p \colon B \to A \mid p \text{ is an operadic fibration in } \dd\Set \text{ between lean $\infty$-operads } A \text{ and } B\}  \]
and the set of generating trivial fibrations is 
\[\mathcal{Q} = \{p \colon B \to A \mid p \text{ is a trivial fibration in } \dd\Set \text{ between lean $\infty$-operads } A \text{ and } B\}. \]
Analogous sets are used in the open and closed case. Since the proofs for open, closed and general dendroidal profinite sets are identical, we will focus on the general case.

\begin{remark}
Since the cofibrations in the operadic model structure on $\dd\wh\Set$ are the normal monomorphisms, one can replace the set $\mathcal{Q}$ with any set of maps in $\dd\wh\Set$ with the property that a map of dendroidal profinite sets $X \to Y$ is normal if and only if it has the \llp \wrt any map in this set. In particular, by \autoref{remark:NormalMonoIffLLPOtherSets}, we could replace $\mathcal{Q}$ by the set of trivial fibrations between weakly contractible lean $\infty$-operads, or by the set of all trivial fibrations between lean dendroidal sets.
\end{remark}

Before embarking on the proof of this theorem, we deduce the following formal consequence.

\begin{proposition}\label{prop:PullbackPowerProperty}
Let $X \cofarrow Y$ be a cofibration and $L \fibarrow K$ a fibration in the model structure for profinite $\infty$-operads on $\dd\wh\Set$. Then the map
\[\sHom(Y,L) \to \sHom(X,L) \times_{\sHom(X,K)} \sHom(Y,K) \]
is a categorical fibration, which is a trivial fibration if either $L \fibarrow K$ or $X \cofarrow Y$ is trivial.
\end{proposition}

\begin{proof}
Since the model structure for profinite $\infty$-operads is fibrantly generated, it suffices to prove this proposition in the case where $L \fibarrow K$ is a generating (trivial) fibration. Let such a generating (trivial) fibration $B \fibarrow A$ be given. In the case that $X \cofarrow Y$ is a trivial cofibration or $B \fibarrow A$ a generating trivial fibration, the proposition follows from \autoref{lemma:PullbackPowerIsJFibration}. The case where $X \cofarrow Y$ is a cofibration and $B \to A$ a generating fibration is equivalent to showing that for any trivial cofibration $M \trivcofarrow N$ in the Joyal model structure on $\s\Set$ and any generating fibration $B \fibarrow A$ of $\dd\wh\Set$, the map
\[B^N \to B^M \times_{A^M} A^N \]
is a trivial fibration in $\dd\wh\Set$. We already know that it is a fibration by one of the other two cases. It is a general fact about model categories that a fibration is a trivial fibration if and only if it has the \rlp \wrt cofibrations between cofibrant objects (see e.g. \cite[Lemma 8.43]{HeutsMoerdijk2020Trees}). It follows by adjunction from \autoref{cor:PullbackPowerBetweenNormal} that $B^N \to B^M \times_{A^M} A^N$ has the \rlp \wrt normal monomorphisms between normal dendroidal profinite sets, so we conclude that it is a trivial fibration.
\end{proof}

The proof of \autoref{theorem:ProfiniteOperadicModelStructure} will use the following lemmas.

\begin{lemma}\label{lemma:TrivialFibrationStableUnderPullback}
For any trivial fibration between lean dendroidal sets $B \trivfibarrow A$ and any map of dendroidal profinite sets $X \to A$, the pullback $X \times_A B \to X$ is a weak equivalence in $\dd\wh\Set$. In particular, any pullback in $\dd\wh\Set$ of a map in the set $\mathcal{Q}$ is a weak equivalence in $\dd\wh\Set$.
\end{lemma}

\begin{proof}
First note that any trivial fibration between lean dendroidal sets is a weak equivalence in $\dd\wh\Set$. Write $X = \{X_i\}_{i \in I}$ with $X_i$ lean and where $I$ is a codirected poset. Since $A$ is lean, the map $X \to A$ factors through $X_i$ for some $i \in I$. Since $I_{\leq i} \subset I$ is cofinal, we may assume without loss of generality that $X \to A$ factors through $X_i$ for every $i \in I$. Then $X \times_A B \to X$ is the inverse limit of the maps $X_i \times_A B \to X_i$. These are all trivial fibrations between lean dendroidal sets, hence weak equivalences. The result now follows from \autoref{prop:WeakEquivalencesStableUnderCofilLimits}.
\end{proof}

\begin{lemma}\label{lemma:NormalMonoOperadicEquivalenceHasLLPwrtFibrations}
Let $X \cofarrow Y$ be a normal monomorphism in $\dd\wh\Set$. If $X \cofarrow Y$ is a weak equivalence, then it has the \llp \wrt operadic fibrations between lean $\infty$-operads, i.e. the maps in $\mathcal{P}$.
\end{lemma}

\begin{proof}
Let $B \fibarrow A$ be a map in $\mathcal{P}$. We will first show that $X \times E \cofarrow Y \times E$ has the \llp \wrt $B \fibarrow A$ and then use this to show that this also holds for $X \cofarrow Y$. Consider the diagram
\[\begin{tikzcd}[column sep = tiny]
\sHom(Y \times E, B) \arrow[drr, bend left = 13] \arrow[ddr, bend right, "\sim" rot320, two heads] \arrow[dr, two heads, "p"] & & \\
& \sHom(X \times E, B) \times_{\sHom(X \times E, A)} \sHom(Y \times E, A) \arrow[dr,phantom, very near start, "\lrcorner"] \arrow[r] \arrow[d, two heads] & \sHom(Y \times E, A) \arrow[d, "\sim" rot90, two heads] \\
& \sHom(X \times E, B) \arrow[r] & \sHom(X \times E, A).
\end{tikzcd}\]
We deduce from \autoref{cor:PullbackPowerBetweenNormal} that $\sHom(Y \times E,A) \to \sHom(X \times E, A)$ is a categorical fibration, while it is a categorical equivalence by definition of the weak equivalences in $\dd\wh\Set$. In particular, it is a trivial fibration, hence its pullback is as well. We furthermore see that $\sHom(Y \times E,B) \to \sHom(X \times E, B)$ is a weak equivalence by definition, so the map $p$ is a categorical equivalence by the 2 out of 3 property. Since it is a categorical fibration by \autoref{cor:PullbackPowerBetweenNormal}, we conclude that this map is a trivial fibration and hence surjective on $0$-simplices. In particular, $X \times E \cofarrow Y \times E$ has the \llp \wrt $A \fibarrow B$.

To deduce from this that $X \cofarrow Y$ has the \llp \wrt $A \fibarrow B$, suppose that maps $X \to B$ and $Y \to A$ are given for which the right-hand square in the diagram
\[\begin{tikzcd}
X \times E \ar[r, two heads, "\sim"] \ar[d, tail]       & X \ar[d, tail] \ar[r]     & B \ar[d, two heads] \\
Y \times E \ar[r, two heads, "\sim"'] \ar[urr, dashed]    & Y \ar[r]                  & A
\end{tikzcd}\]
commutes. By what we proved above, the indicated lift exists. By \autoref{lem:NormalMonoHasIncreasinglyNormalRepresentation}, we may assume without loss of generality that $X \to Y$ has an increasingly normal representation $\{X_i \to Y_i\}_{i \in I}$. Since $A$ and $B$ are lean, there exists an $i$ such that $X \to B$ and $Y \to A$ factor through $X_i$ and $Y_i$, respectively, the lift $Y \times E \to B$ factors through $Y_i \times E$, and the diagram
\[\begin{tikzcd}
X_i \times E \ar[r, two heads, "\sim"] \ar[d, tail]       & X_i \ar[d, tail] \ar[r]     & B \ar[d, two heads] \\
Y_i \times E \ar[r, two heads, "\sim"'] \ar[urr, dashed]    & Y_i \ar[r]                  & A
\end{tikzcd}\]
commutes. If we can construct a lift in the right-hand square, then it follows that the desired lift $Y \to B$ exists. By choosing $i$ small enough, we may assume that $X_i \to Y_i$ is $n$-normal, where $n$ is some number such that $A$ and $B$ are $n$-coskeletal. Finding a lift in the right-hand square of the above diagram is then equivalent to constructing a lift in the right-hand square of the following diagram:
\begin{equation}\label{diagram:skeletalLiftingPropblem}
\begin{tikzcd}
\sk_n (X_i) \times E \ar[rr, two heads, "\sim"] \ar[dd, tail]     & & \sk_n(X_i) \ar[dd, tail] \ar[dl, tail] \ar[rr]     & & B \ar[dd, two heads] \\
 & |[yshift=-1.2em, xshift=1.2em, overlay]|  P \ar[dr] \ar[urrr,dashed,"l"'] \arrow[ul,phantom, very near start, "\ulcorner"] & & \\
\sk_n(Y_i) \times E \ar[rr, two heads, "\sim"'] \ar[ur, "\sim" rot20]    & & \sk_n(Y_i) \ar[rr]  & & A.
\end{tikzcd}
\end{equation}
Here $\sk_n(X_i) \to \sk_n(Y_i)$ is normal by \autoref{lem:nSkeletonOfNNormalMapIsNormal}, and $P$ is defined as a pushout. The universal property of the pushout provides us with the dashed lift $l \colon P \to B$. Since the operadic model structure on $\dd\Set$ is left proper, the map $\sk_n(Y_i) \times E \to P$ is an operadic weak equivalence, hence $P \to \sk_n(Y_i)$ is an operadic weak equivalence by the 2 out of 3 property. Factor this map as a trivial cofibration followed by a trivial fibration: $P \trivcofarrow Z \trivfibarrow \sk_n(Y_i)$. The operadic model structure on $\dd\Set$ provides us with two lifts
\[
\begin{tikzcd}
\sk_n(X_i) \ar[d, tail] \ar[r] & Z \ar[d, two heads, "\sim" rot90] \\
\sk_n(Y_i) \ar[ur,dashed] \ar[r,"="'] & \sk_n(Y_i)
\end{tikzcd}
\qquad \text{and} \qquad\quad
\begin{tikzcd}
P \ar[d,tail,"\sim" rot90] \ar[r, "l"] & B \ar[d,two heads] \\
Z \ar[ur, dashed] \ar[r] & A.
\end{tikzcd}
\]
Composing these two lifts provides a map $\sk_n(Y_i) \to B$ that is a lift in the right-hand square of diagram \eqref{diagram:skeletalLiftingPropblem}, so we conclude that the desired lift $Y \to B$ exists.
\end{proof}

We are now ready to prove the existence of the model structures.

\begin{proof}[Proof of \autoref{theorem:ProfiniteOperadicModelStructure}]
We check the hypotheses of the dual of \cite[Theorem 11.3.1]{Hirschhorn2003Model}. The weak equivalences in $\dd\wh\Set$ satisfy the two out of three property and are closed under retracts since this holds for the categorical equivalences in $\s\Set$.

\par (1) Since lean dendroidal sets are cocompact in $\dd\wh\Set$, the sets $\mathcal{P}$ and $\mathcal{Q}$ permit the cosmall object argument.

\par (2) As mentioned above, weak equivalences in $\dd\wh\Set$ are closed under retracts. Moreover, by \autoref{prop:WeakEquivalencesStableUnderCofilLimits} the weak equivalences in $\dd\wh\Set$ are closed under cofiltered limits, so a transfinite precomposition of weak equivalences is again an weak equivalences. Therefore, in order to show that any map in the cosaturation of $\mathcal{Q}$ is a weak equivalence, it suffices to show that any pullback in $\dd\wh\Set$ of a map in $\mathcal{Q}$ is a weak equivalence. This follows from \autoref{lemma:TrivialFibrationStableUnderPullback}.

\par (3) Any map $X \to Y$ that has the \llp \wrt maps in $\mathcal{P}$ does so with respect to maps in $\mathcal{Q}$, since $\mathcal{Q} \subset \mathcal{P}$. Furthermore, such a map is a weak equivalence by \autoref{prop:NormalMonoOperadicEquivalenceIfLLPwrtFibrations}.

\par (4) We need to show that if $X \to Y$ has the \llp \wrt maps in $\mathcal{Q}$ and is a weak equivalence, then $X \to Y$ has the \llp \wrt all maps in $\mathcal{P}$. This follows from \autoref{lemma:NormalMonoOperadicEquivalenceHasLLPwrtFibrations}, noting that $X \to Y$ is normal by \autoref{prop:NormalMonoIffLLPwrtLeanTrivFibs}.

Finally, for left properness, let a pushout square of the form
\[\begin{tikzcd}
X \ar[r, "\sim"] \ar[d, tail] & Y \ar[d] \\
Z \ar[r] & W \ar[ul, phantom, very near start, "\ulcorner"]
\end{tikzcd}\]
be given. Note that the cartesian product of a finite set $A$ and profinite set $K$ is simply the coproduct $\coprod_{a \in A} K$, hence pushouts of profinite sets are preserved by cartesian products with finite sets. Since $E$ is degreewise finite and colimits in $\dd\wh\Set$ are computed degreewise, it follows that $W \times E \cong Y \times E \sqcup_{X \times E} Z \times E$. Since $X \times E \wearrow Y \times E$ is a weak equivalence between cofibrant objects, the pushout $Z \times E \to W \times E$ must also be a weak equivalence. We conclude that $Z \to W$ is a weak equivalence.
\end{proof}

\begin{remark}\label{remark:CovariantPicard}
	A careful look at the proof of \autoref{theorem:ProfiniteOperadicModelStructure} shows that it only uses a few properties of the operadic model structure, which are also enjoyed by many other model structures on the category of (open or closed) dendroidal sets. For example, one can use the same method to obtain profinite versions of the covariant model structure of \cite[\S 9.6]{HeutsMoerdijk2020Trees} and the Picard model structure of \cite{BasicNikolaus2014PicardModelStructure} (called the \emph{stable} model structure there).
\end{remark}

Let us observe some immediate consequences of the theorem. Recall from \autoref{ssec:ProCategories} that the inclusion $\LdSet \hookrightarrow \dd\Set$ extends to a limit-preserving functor $U \colon \dd\wh\Set \to \dd\Set$ whose left adjoint is the profinite completion functor, and similar for open and closed dendroidal sets.

\begin{corollary}
The profinite completion functors $\dd\Set \to \dd\wh\Set$, $\od\Set \to \od\wh\Set$ and $\cd\Set \to \cd\wh\Set$ are left Quillen with respect to the operadic model structures on their domain and the model structures for profinite (open/closed) $\infty$-operads on their target.
\end{corollary}

\begin{proof}
The right adjoint makes the diagram
\[\begin{tikzcd}
\LdSet \ar[d, hook] \ar[dr, hook] & \\
\dd\wh\Set \ar[r] & \dd\Set \\
\end{tikzcd}\]
commute, so it must send generating (trivial) fibrations to (trivial) fibrations in $\dd\Set$, and similar for the open and closed case.
\end{proof}

One can prove the following fact about the fibrant objects in the model structure for profinite $\infty$-operads.

\begin{proposition}
Any fibrant object in the model structure for profinite (open/closed) $\infty$-operads is an inverse limit of lean (open/closed) $\infty$-operads.
\end{proposition}

\begin{proof}
This proof is dual to that of Lemma 7.11 in \cite{BlomMoerdijk2020SimplicialProV1}.
\end{proof}

\begin{remark}\label{remark:AssociativeOperadNotProfinite}
	One can deduce from the previous proposition that the nerve of the unital associative operad $\mathcal{A}ss$ is not fibrant in the model structure for profinite $\infty$-operads on $\dd\wh\Set$. If it were, then it could be written as an inverse limit $\lim_i A_i$ of lean $\infty$-operads. We obtain ordinary operads $\mathcal{P}_i$ by applying $\pi_0$ to all spaces of operations of these lean $\infty$-operads $A_i$, analogous to how one obtains the homotopy category of an $\infty$-category. Then $\mathcal{A}ss = \lim_i \mathcal{P}_i$. Since $\mathcal{A}ss$ has one colour, the projection maps $\mathcal{A}ss \to \mathcal{P}_i$ all land in a single colour of $\mathcal{P}_i$. By discarding all other colours in $\mathcal{P}_i$, we may assume without loss of generality that the operads $\mathcal{P}_i$ are uncoloured. For any $i$, the fact that $A_i$ is lean implies that there exists an $n_i$ such that $\mathcal{P}_i(n) = *$ for $n \geq n_i$. Using the nullary operation $* \in \mathcal{A}ss(0)$ and the $\Sigma_m$-equivariantness of $\mathcal{A}ss(m) \to \mathcal{P}_i(m)$, this implies that for each $i$, the projection $\mathcal{A}ss \to \mathcal{P}_i$ must factor through the terminal object (i.e. the commutative operad). In particular, if $\mathcal{A}ss$ is isomorphic to the inverse limit $\lim_i \mathcal{P}_i$, then the identity map of $\mathcal{A}ss$ would factor through the terminal object, which is of course not possible.
\end{remark}

Let us call a map $\wt X \to X$ in $\dd\wh\Set$ a \emph{normalization of $X$} if $\wt X$ is normal and $\wt X \to X$ is a weak equivalence. For any object $X$ in $\dd\wh\Set$, the trivial fibration $X \times E \to X$ used to define the weak equivalences is such a normalization. It follows from the existence of the model structure that the same weak equivalences can be characterized by an arbitrary (functorial) normalization:

\begin{proposition}\label{prop:AlternativeNormalizations}
	Consider a map $f \colon X \to Y$ in $\dd\wh\Set$ and normalizations which fit into a commutative square
	\[\begin{tikzcd}
		\wt X \ar[r,"\wt f"] \ar[d, "\sim" rot90, "p"] & \wt Y \ar[d, "\sim" rot90,"q"] \\
		X \ar[r,"f"] & Y.
	\end{tikzcd}\]
	Then $f \colon X \to Y$ is a weak equivalence if and only if for each lean $\infty$-operad $A$, the map $\sHom(\wt Y, A) \to \sHom(\wt X, A)$ is an equivalence of $\infty$-categories.
\end{proposition}

\begin{proof}
	By Brown's lemma and \autoref{prop:PullbackPowerProperty}, the map $\sHom(\wt Y, A) \to \sHom(\wt X, A)$ is a weak equivalence if and only if $\sHom(\wt Y \times E, A) \to \sHom(\wt X \times E, A)$ is, so the result follows.
\end{proof}

Recall from \autoref{ssec:dSetsAndSpaces} that the adjunctions
\[\begin{tikzcd}
	\dd\Set \arrow[r, "o^*"', bend right, shift left = 1, ""{name=B, above}] & \od\Set, \arrow[l, "o_!"', bend right, shift left = 1, ""{name=A, below}] \ar[phantom, from=A, to=B, "\dashv" rotate=-90, no line]
\end{tikzcd}\quad
\begin{tikzcd}
	\dd\Set \arrow[r, "u^*"', bend right, shift left = 1, ""{name=B, above}] & \cd\Set \arrow[l, "u_!"', bend right, shift left = 1,""{name=A, below}] \ar[phantom, from=A, to=B, "\dashv" rotate=-90, no line]
\end{tikzcd} \quad \text{and} \quad
\begin{tikzcd}
	\od\Set \arrow[r, "h_*"', bend right, shift left = 1, ""{name=B, above}] & \cd\Set, \arrow[l, "h^*"', bend right, shift left = 1, ""{name=A, below}] \ar[phantom, from=A, to=B, "\dashv" rotate=-90, no line]
\end{tikzcd}\]
induced by $o \colon \Omegao \to \Omega$, $u \colon \Omegacl \to \Omega$ and $h \colon \Omegao \to \Omegacl$ are Quillen pairs. Since $\wh\Set$ is (co)complete, we obtain analogous adjunctions defined in terms Kan extensions in the profinite case. These can also be shown to be Quillen pairs.

\begin{proposition}
	The functors
	\[o^* \colon \dd\wh\Set \to \od\wh\Set, \quad u^* \colon \dd\wh\Set \to \cd\wh\Set \quad \text{and} \quad h_* \colon \od\wh\Set \to \cd\wh\Set \]
	are right Quillen functors with respect to the model structures for profinite (open/closed) $\infty$-operads
\end{proposition}

\begin{proof}
	By Remarks \ref{remark:PushforwardPreservesLean} and \ref{remark:PullbackOpenPreservesLean}, the ordinary versions of the functors $o^*$, $u^*$ and $h_*$ preserve lean objects, hence they agree with their profinite analogues on lean (open) dendroidal sets. In particular, these functors preserve the generating (trivial) fibrations of the model structures on $\dd\wh\Set$ and $\od\wh\Set$.
\end{proof}

Recall from \autoref{ssec:dSetsAndSpaces} that the functor $o_!$ can be used to identify $\od\Set$ with the over-category $\dd\Set/O$, and that under this identification their model structures agree. Since $O$ is a degreewise finite dendroidal set, it can also be seen as an object of $\dd\wh\Set$, and one similarly obtains that $\od\wh\Set \simeq \dd\wh\Set/O$.

\begin{proposition}\label{prop:OpenDendroidalProfiniteSetsIsSlice}
	The operadic model structure on $\od\wh\Set$ agrees with the model structure on $\dd\wh\Set/O$ obtained by slicing the model structure for profinite $\infty$-operads.
\end{proposition}

\begin{proof}
	Since the equivalence of categories $\od\wh\Set \simeq \dd\wh\Set/O$ is established by the functor $o_! \colon \od \wh\Set \to \dd\wh\Set$, it suffices to show that a map $f$ in $\od\wh\Set$ is a (co)fibration if and only if $o_! f$ is. For cofibrations, this follows from the ordinary case since a map in $\od\wh\Set$ or $\dd\wh\Set$ is normal if and only if the underlying map is normal in $\dd\Set$ or $\od\Set$. For fibrations, note that $o^* \colon \dd\wh\Set/O \to \od\wh\Set$ is right Quillen and inverse to $o_!$, so it suffices to show that $o_!$ preserves fibrations. By \autoref{remark:PullbackOpenPreservesLean} the ordinary version of $o_!$ preserves lean objects, hence it agrees with the profinite version of $o_!$ on lean objects. Since the ordinary version of $o_!$ preserves fibrations, this implies that $o_! \colon \od\wh\Set \to \dd\wh\Set/O$ preserves generating fibrations, hence all fibrations.
\end{proof}

If one slices over $\eta$ instead, then one obtains the equivalence $\dd\Set/\eta \simeq \s\Set$ and the analogous equivalence in the profinite case. For ordinary dendroidal sets, it can be shown that the operadic model structure restricts to the Joyal model structure under this equivalence. By slicing the operadic model structure on $\dd\wh\Set$ over $\eta$, one similarly obtains the profinite Joyal model structure of \cite[Example 5.7]{BlomMoerdijk2020SimplicialProV1}.

\begin{proposition}\label{prop:SlicingProfiniteOperadsOverEta}
	The profinite Joyal model structure on $\s\wh\Set$ agrees with the model structure on $\dd\wh\Set/\eta$ obtained by slicing the model structure for profinite $\infty$-operads.
\end{proposition}

\begin{proof}
	By \autoref{prop:OpenDendroidalProfiniteSetsIsSlice}, it suffices to show that the model structure on $\od\wh\Set/\eta$ agrees with the profinite Joyal model structure on $\s\wh\Set$ under the equivalence $i^* \colon \od\wh\Set/\eta \to \s\wh\Set$. Note that $i^*$ is an equivalence of $\s\Set$-enriched categories, where the enrichment of $\od\wh\Set/\eta$ comes from the simplicial hom of $\od\wh\Set$.
	
	It is clear that the cofibrations agree, since these are the monomorphisms in both cases. By \autoref{remark:UnderlyingInftyCatPreservesLean}, $i^*$ sends lean objects in $\od\wh\Set$ to lean objects in $\s\wh\Set$. Since the ordinary version of $i^*$ preserves fibrations, we see that $i^*$ must send all generating fibrations of $\od\wh\Set$ to (generating) fibrations in $\s\wh\Set$ by definition of the fibrations in the profinite Joyal model structure (cf. \cite[Example 5.7]{BlomMoerdijk2020SimplicialProV1}). This implies that its inverse $i_!$ preserves trivial cofibrations and hence weak equivalences since all objects are cofibrant.
	
	Conversely, to show that $i^*$ preserves weak equivalences, let $X \to Y$ be a weak equivalence in $\od\wh\Set/\eta$. We need to show that $\sHom(i^*Y, Z) \to \sHom(i^*X, Z)$ is a weak equivalence of $\infty$-categories for any lean $\infty$-category $Z$. Let $m$ be such that $Z$ is $m$-coskeletal. Now note that $i^*\cosk_{2m}(i_! Z) = Z$, hence showing that the above map is an equivalence of $\infty$-categories is equivalent to showing that $\sHom(Y,\cosk_{2m}(i_! Z)) \to \sHom(X,\cosk_{2m}(i_! Z))$ is. Since $\cosk_{2m}(i_! Z)$ is a lean open dendroidal set, this follows if we can show that it is fibrant in $\od\Set$.
	
	By adjunction, this holds if and only if $i_! Z$ has the \rlp \wrt the inclusion $\sk_{2m} \Lambda^e[T] \to \sk_{2m} \Omega[T]$ for every tree $T$. If $|T| > 2m+2$, then this map is an isomorphism so there is nothing to prove. If $|T| \leq 2m$, then this is simply the horn $\Lambda^e[T] \cofarrow \Omega[T]$, so it holds since $i_! Z$ is an open $\infty$-operad. In the cases where $|T| = 2m+1$ or $|T|= 2m+2$ one has $\sk_{2m} \Lambda^e[T] = \Lambda^e[T]$, using that $T$ is open in the case that $|T| = 2m+2$. In particular, given a map $\sk_{2m} \Lambda^e[T] \to i_! Z$, one can extend this to a map $\Omega[T] \to i_!Z$ and then restrict it to a map $\sk_{2m}\Omega[T] \to i_! Z$ to construct the desired lift. We conclude that $\cosk_{2m} i_! Z$ is a lean $\infty$-operad.
\end{proof}

\printbibliography

\noindent{\sc Matematiska institutionen, Stockholms universitet, 106 91 Stockholm, Sweden.}\\
\noindent{\emph{E-mail:} \href{mailto:blom@math.su.se}{\nolinkurl{blom@math.su.se}}}\vspace{2ex}

\noindent{\sc Mathematisch Instituut, Universiteit Utrecht, Postbus 80010, 3508 TA Utrecht, The Netherlands.}\\
\noindent{\emph{E-mail:} \href{mailto:i.moerdijk@uu.nl}{\nolinkurl{i.moerdijk@uu.nl}}}

\end{document}